\documentclass[11pt]{article}
\usepackage{latexsym}
\usepackage{amsmath}
\usepackage{cite}
\usepackage{amsthm,color,bm}
\usepackage{amssymb}
\usepackage{mathrsfs}
\usepackage{indentfirst}
\usepackage{times,cases}
\usepackage{dsfont}

\usepackage{lscape}
\newtheorem{theorem}{\noindent Theorem}[section]
\newtheorem{proposition}[theorem]{\noindent Proposition}
\newtheorem{definition}[theorem]{\noindent Definition}
\newtheorem{lemma}[theorem]{\noindent Lemma}
\newtheorem{remark}[theorem]{\noindent Remark}
\newtheorem{corollary}[theorem]{\noindent Corollary}
\numberwithin{figure}{section}
\numberwithin{equation}{section}

\renewcommand{\theequation}{\thesection.\arabic{equation}}

\DeclareMathOperator*{\esssup}{ess\,sup}
\DeclareMathOperator*{\LIM}{LIM}

\newcommand{\cA}{{\mathcal A}}
\newcommand{\cB}{{\mathcal B}}
\newcommand{\cC}{{\mathcal C}}
\newcommand{\cD}{{\mathcal D}}

\newcommand{\cF}{{\mathcal F}}

\newcommand{\cI}{{\mathcal I}}

\newcommand{\cM}{{\mathcal M}}

\newcommand{\cP}{{\mathcal P}}
\newcommand{\cQ}{{\mathcal Q}}
\newcommand{\cR}{{\mathcal R}}

\newcommand{\cT}{{\mathcal T}}
\newcommand{\cU}{{\mathcal U}}

\newcommand{\sM}{{\mathscr M}}

\def\Z{\mathbb{Z}}
\def\R{\mathbb{R}}

\def\bP{\mathbb{P}}

\def\N{\mathbb{N}}
\def\1{\mathds{1}}
\newcommand{\di}{\mathrm{d}}
\newcommand{\me}{\mathrm{e}}

\def\disp{\displaystyle}

\def\bc{\begin{center}}
\def\ec{\end{center}}
\def\be{\begin{equation}}
\def\ee{\end{equation}}
\def\bea{\begin{eqnarray}}
\def\eea{\end{eqnarray}}
\def\ba{\begin{array}}
\def\ea{\end{array}}
\def\benu{\begin{enumerate}}
\def\eenu{\end{enumerate}}
\def\bt{\begin{theorem}}
\def\et{\end{theorem}}
\def\bl{\begin{lemma}}
\def\el{\end{lemma}}
\def\bco{\begin{corollary}}
\def\eco{\end{corollary}}
\def\bn{\begin{numcases}}
\def\en{\end{numcases}}
\def\br{\begin{remark}}
\def\er{\end{remark}}
\def\bd{\begin{definition}}
\def\ed{\end{definition}}
\def\bp{\begin{proposition}}
\def\ep{\end{proposition}}
\def\bo{\begin{proof}}
\def\eo{\end{proof}}
\def\bx{\begin{example}}
\def\ex{\end{example}}
\def\bal{\begin{align}}
\def\eal{\end{align}}

\def\pa{\partial}
\def\al{\alpha}
\def\De{\Delta} \def\de{\delta}

\def\lam{\lambda} 

\def\ve{\varepsilon}
\def\sig{\sigma}
\def\vsig{\varsigma}
\def\vp{\varphi}

\def\w{\omega}\def\W{\Omega}
\def\gam{\gamma}

\def\~{\widetilde}

\def\ol{\overline}

\def\Cup{\bigcup}

\def\ra{\rightarrow}
\def\Ra{\Rightarrow}

\def\8{\infty}
\def\X{\times}

\def\mb{\mbox}
\def\di{{\rm d}}
\def\me{{\rm e}}

\def\suo{\!\!\!}

\def\Hs{\hspace{0.8cm}}
\def\hs{\hspace{0.4cm}}
\def\Vs{\vskip10pt}
\def\vs{\vskip5pt}

\def\({\left(}
\def\){\right)}

\begin{document}


\begin{center}
    {\large \bf Invariant sample measures and sample statistical solutions for\\
    nonautonomous stochastic lattice Cahn-Hilliard equation with nonlinear noise}
\vspace{0.5cm}\\
{Jintao Wang*,\hs Dongdong Zhu,\hs Chunqiu Li}\\\vspace{0.3cm}

{\small Department of Mathematics, Wenzhou University, Wenzhou 325035, China}
\end{center}


\renewcommand{\theequation}{\arabic{section}.\arabic{equation}}
\numberwithin{equation}{section}


\begin{abstract}
We consider a stochastic lattice Cahn-Hilliard equation
with nonautonomous nonlinear noise.
First, we prove the existence of pullback random attractors in $\ell^2$
for the generated nonautonomous random dynamical system.
Then, we construct the time-dependent invariant sample Borel probability measures
based on the pullback random attractor.
Moreover, we develop a general stochastic Liouville type equation
for nonautonomous random dynamical systems and
show that the invariant sample measures obtained satisfy the stochastic Liouville type equation.
At last, we define a new kind of statistical solution --- sample statistical solution
corresponding to the invariant sample measures
and show that each family of invariant sample measures is
a sample statistical solution.
\vs

\noindent\textbf{Keywords:} Stochastic lattice Cahn-Hilliard system; pullback random attractors;
invariant sample measures; stochastic Liouville type theorem; statistical solutions.
\vs

\noindent{\bf AMS Subject Classification 2010:}\, 60H10, 37L60, 37L30, 37L40

\end{abstract}

\vspace{-1 cm}

\footnote[0]{\hspace*{-7.4mm}
$^{*}$ Corresponding author.\\
E-mail address: wangjt@wzu.edu.cn (J.T. Wang); 22451025027@stu.wzu.edu.cn (D.D. Zhu);
lichunqiu@wzu.edu.cn (C.Q. Li).}

\section{Introduction}\label{s1}

This article mainly concerns with long-term dynamical behaviors of
the following nonautonomous stochastic Cahn-Hilliard equation with infinite lattices
and nonautonomous nonlinear noise from the invariant sample measures and
sample statistical solutions:
\begin{align}\label{1.1}
\di u_i(t)+[(A(Au+f(u)))_i+\lam u_i-g_i(t,u_i)]\di t=(k u_i+h_i(t))\di W(t),\hs t>\tau
\end{align}
with the initial datum
\be\label{1.2}u_i(\tau)=u_{\tau i}\hs\mb{and}\hs i\in\Z,\ee
where the linear mapping $A:\ell^2\ra\ell^2$ is defined as
\be\label{1.3}(Au)_i=2u_i-u_{i-1}-u_{i+1},\ee
$f:\ell^2\ra\ell^2$ is a mapping, $\lam>0$, $g(t,u)=(g_i(t,u_i))_{i\in\Z}$
is the nonautonomous nonlinear forcing,
$k\geqslant 0$, $h(t)=(h_i(t))_{i\in\Z}$,
$W$ is a two-sided real-valued Wiener process on the probability space to be determined later
and $u_{\tau}:=(u_{\tau,i})_{i\in\Z}\in\ell^2$.
The lattice equation \eqref{1.1} can be regarded as the spatial discrete form
on infinite lattices of 1-dimensional nonautonomous stochastic Cahn-Hilliard equation
\be\label{1.4}\frac{\pa u}{\pa t}+(u_{xx}+f(u))_{xx}+\lam u
-g(t,u)=(ku+h(t))\frac{\di W(t)}{\di t},\hs x\in\R.\ee

Cahn-Hilliard equation arose originally in materials science
to describe rapid decomposition phenomena of
a binary alloy system for the 1D case by Cahn and Hilliard \cite{CH58}.
It also pertains to the Model B class in Hohenberg and Halperin's classification \cite{HH77}.
This is a standard model for phase transition with conserved quantities and
can be applied to phase transition in liquid crystals \cite{CCCG},
segregation of granular mixtures in a rotating drum \cite{PH},
formation of sand ripples \cite{SMM99,SW99}, coalescence phenomena \cite{V-G04} and
cell diffusion, proliferation and adhesion in biology \cite{CMZ14,D16,KS08,M13}.
Moreover, Cahn-Hilliard equation is also a partial differential equation to which a conservative
noise is added to account for thermal fluctuations \cite{C70}.
As further applications, the immiscible binary fluids can be simulated
by (singular) Cahn-Hilliard equation (\cite{Ch18,FM23,LT98}),
which hence can be used to study microconcentrations
and microstructure evolution of the solder alloy (\cite{FM23,USG04}).

In recent several decades, Cahn-Hilliard equations have been studied from many distinct views.
The posedness problems were discussed for different sorts of Cahn-Hilliard equations
in \cite{CM95,EZ86,GMS09,L22,S21,WLu19}.
Goh and his team studied Hopf bifurcations of Cahn-Hilliard equation in a series of works
including \cite{GS15}.
The long-term asymptotic behavior fascinated a lot more attention these years;
see \cite{CG16,CMZ11,D16,GMS10,GSZ09,JWLD11,LZ98,M13,Sc07,Se07}.
For autonomous Cahn-Hilliard with stochastic noise, invariant measure was also considered
(\cite{JSW14,QAH}).
Moreover, there have been some works using the lattice Boltzmann models to study
the Cahn-Hilliard equation for numerical simulation in \cite{LSW12,LSGC14,WCSL16,ZLWC24},
where, however, the framework is far different from
the discretized system with infinite lattices presented in \cite{HSZ11,LLW22,LLW23,WJ24}.
By comparison, the investigations of lattice dynamical systems
and invariant (sample) measures for nonautonomous Cahn-Hilliard systems,
which will be constructed in this article,
have been almost blank up to our knowledge.

Lattice dynamical systems are developed by spatial discretizations of
general partial differential equations,
where a countable system of ordinary differential equations is constructed
by replacing the spatial derivatives by differences (\cite{D77,WJ24,WWHK}).
The applications of lattice systems have spread in many important areas;
see \cite{CMV14,HSZ11,LLW22,LLW23,WJ24} and their references.
Invariant measures and statistical solutions (satisfying a Liouville type equation)
perform a significant role in the research of turbulence (see Foias et al. \cite{FMRT01}),
which is a very important observation target in fluid dynamics,
and invariant measures have been studied for many classical model for fluid dynamics
(see \cite{LLW22,LLW23,WJ24,WZ23,WZC20,WZZ21,WWL22}).

Invariant sample measures for nonautonomous random dynamical systems (NRDS's for short)
were recently developed in \cite{CY23,WJ24,ZWC22}
from time-dependent invariant measure for deterministic systems.
This sort of invariant measure is a ``weaker" version of the usual one for
random dynamical systems,
since the "sample" here represents the reliance of the invariant measures on each fixed $\w\in\W$,
while the usual one is in the sense of expectation in $\W$ (\cite{WLYJ21,WZL23}).
In spite of the weak sense for NRDS's, the invariant sample measures
still deserve special attention, in light of the difficulty to extend the classical methods for
random dynamical systems to NRDS's (see \cite{WZ23}) to obtain
the time-dependent invariant measures in the usual sense.
As an important property, invariant sample measures was proved by Chen and Yang in \cite{CY23}
to satisfy a stochastic Liouville type equation, where they introduced a class of cylindrical
test functions and used It\^o's formula.
But until now, there has not yet been any appropriate statistical solution defined
with respect to the invariant sample measures for a given NRDS.

In this article, we intend to develop the \emph{sample statistical solutions}
by means of the invariant sample measures for the lattice system generated
by nonautonomous stochastic Cahn-Hilliard equation \eqref{1.1} with nonlinear noise.
In order to achieve a wider range of applications,
we select the pivotal points of the cylindrical test functions introduced in \cite{CY23},
and further sum up four key conditions that should be satisfied by a test function.
Then we put up a more general definition of the class of test functions as follows.

Let $H$ be a sparable Hilbert space with inner product $(\cdot,\cdot)$ and norm $\|\cdot\|$
and $V$ a Banach space with $V^*$ its dual space and dense embeddings $V\subset H\subset V^*$.
Denote the dual pairing between $v^*\in V^*$ and $v\in V$ by $\langle v^*,v\rangle$.
We consider the following stochastic differential equation
\be\label{1.5}\di u=\tilde{F}(t,u)\di t+\tilde{G}(t,u)\di W,\ee
in $H$, where $\tilde{F}:\R\X V\ra V^*$ is the drift term,
$\tilde{G}:\R\X H\ra H$ is the diffusion coefficient
and $W$ is given in Subsection \ref{ss2.2}.
Assume that for each initial time $\tau\in\R$, $\bP$-a.s. $\w\in\W$ and initial datum $\phi\in H$,
\eqref{1.5} has a global solution $u(\cdot,\tau,\w,\phi)\in\cC^1([\tau,\8);H)$,
where the ergodic metric dynamical system $(\W,\cF,\bP,\{\theta_t\}_{t\in\R})$ is defined
in Subsection \ref{ss2.2}.

\bd\label{de1.1} Let $\cT$ denote the \textbf{class of test functions} $\Psi:H\ra\R$
for \eqref{1.5}
to be bounded on bounded subsets of $H$ such that the following conditions hold:
\benu\item[(1)] for each $u\in V$, the Fr\'echet derivative $\Psi'(u)$
taken in $H$ along $V$ exists.
More precisely, for each $u\in V$, there exists an element in $H$ denoted by $\Psi'(u)$ such that
$$\frac{|\Psi(u+v)-\Psi(v)-(\Psi'(u),v)|}{\|v\|}\ra0\hs\mb{as }\|v\|\ra0,\;v\in V;$$
\item[(2)] for each $u\in V$, $\Psi'(u)\in V$ and the mapping $u\mapsto\Psi'(u)$ is continuous
and bounded as a functional from $V$ into $V$;
\item[(3)] for each $u\in V$, the second-order Fr\'echet derivative $\Psi''(u)$ is
a bounded bilinear operator from $H\X H$ to $\R$ and the mapping $u\mapsto\Psi''(u)$ is continuous
and bounded from $V$ into $L(H\X H,\R)$;
\item[(4)] for every global solution $u(t)$ of \eqref{1.5}, the following
It\^o's formula holds for all $t\geqslant\sig$,
\begin{align}\Psi(u(t))-\Psi(u(\sig))=&\int_\sig^t\langle\tilde{F}(\vsig,u(\vsig)),
\Psi'(u(\vsig))\rangle\di\vsig
+\int_{\sig}^t\(\tilde{G}(\vsig,u(\vsig)),\Psi'(u(\vsig))\)\di W(\vsig)\notag\\
&+\frac12\int_{\sig}^t\Psi''(u(\vsig))\(\tilde{G}(\vsig,u(\vsig)),\tilde{G}(\vsig,u(\vsig))\)
\di\vsig.
\label{1.6}\end{align}
\eenu\ed

Bringing this class $\cT$ of test functions into our framework,
we define a general \emph{stochastic Liouville type equation} for the stochastic
differential equation \eqref{1.5} as, for fixed $(\tau,\w)\in\R\X\W$ and all $t\geqslant s$,
\begin{align}&\int_{H}\Psi(u)\mu_{t,\theta_{t-\tau}\w}(\di u)
-\int_{H}\Psi(u)\mu_{s,\theta_{s-\tau}\w}(\di u)\notag\\
=&\int_s^t\int_H\langle\tilde{F}(\vsig,u),\Psi'(u)\rangle
\mu_{\vsig,\theta_{\vsig-\tau}\w}(\di u)\di\vsig
+\int_s^t\int_H(\tilde{G}(\vsig,u),\Psi'(u))\mu_{\vsig,\theta_{\vsig-\tau}\w}(\di u)
\di\tilde{W}(\vsig)\notag\\
&+\frac12\int_s^t\int_H\Psi''(u)(\tilde{G}(\vsig,u),\tilde{G}(\vsig,u))
\mu_{\vsig,\theta_{\vsig-\tau}\w}(\di u)\di\vsig,\label{1.7}
\end{align}
where $\tilde{W}(\cdot)=W(-\tau+\cdot)-W(-\tau)$, $\Psi\in\cT$ and
$\{\mu_{\tau,\w}\}_{(\tau,\w)\in\R\X\W}$ is a family of invariant sample measures.
Ulteriorly, based on the stochastic Liouville type equation,
we are enabled to develop the statistical solutions into stochastic cases
and define the sample statistical solutions.
We accordingly prove that each family of invariant sample measures constructed
in $\ell^2$ for \eqref{1.1} is a sample statistical solution in our new framework.

In order to construct the invariant sample measures and sample stochastic solutions
for nonautonomous stochastic Cahn-Hilliard equation with nonlinear noise,
we first prove existence of pullback random attractors.
In this procedure, we are necessary to transform the stochastic
differential equation to a random one by using the stationary stochastic solution \eqref{2.5n}
of the Ornstein-Uhlenbeck (OU for short) equation.
However, the transformation methods are different for the cases when $k>0$ and
when $k=0$ (nonautonomous additive noise),
and hence we have to employ distinct estimations to obtain corresponding results
in two sections, respectively.
In both cases, for well-posedness of \eqref{1.1}, we only require that
the partial derivatives $\frac{\pa g_i(t,u_i)}{u_i}$ of the nonautonomous nonlinearity
$g_i(t,u_i)$ have a common upper bound in $L^{\8}_{\rm loc}(\R)$ for all $i\in\Z$,
which means that $g(t,u)$ may not be locally Lipschitzian.
But due to affect of $A$ upon the nonlinearity $f$, $f$ is still necessary to
be locally Lipschitzian.

After constructing the invariant sample measures $\{\mu_{\tau,\w}\}_{(\tau,\w)\in\R\X\W}$,
when checking that $\mu_{\tau,\w}$ satisfy the stochastic Liouville type equation,
one needs to guarantee $\mu_{\vsig,\theta_{\vsig-\tau}\w}$-integrability of the integrands
and further $L^1$-integrability and stochastic integrability of corresponding
$\mu_{\vsig,\theta_{\vsig-\tau}\w}$-integrals in \eqref{1.7}.
It is often a intractable problem how to check the above integrability,
which can not be obtained directly by the settings in Definition \ref{de1.1} after all.
However, in our present problem, thanks to the necessary setting that $V$ is indeed
the space $H$ itself, we can thus overcome this problem.
Consequently, based on the conclusions proved above,
we can check each family of invariant sample measures to be a sample
statistical solution by a standard method.

The rest of this article is organized as follows.
In Section 2, we present the basic settings of this article and the related concepts
and theorems, including some explanation of Definition \ref{de1.1} and
stochastic Liouville type equation \eqref{1.7}.
In Section 3, we study \eqref{1.1} for the case when $k>0$, including the well-posedness,
existence of pullback random attractors, invariant sample measures, stochastic Liouville
type equation and sample statistical solutions.
In Section 4, we consider the similar results for \eqref{1.1} for the case when $k=0$.
Section 5 is the summary and remarks of this article.

\section{Preliminaries}
\subsection{Basic settings}
We introduce the basic settings of spaces and operators
that will be used frequently in this article.

For metric spaces $X$ and $Y$,
we conventionally denote $\cC(X,Y)$ ($\cC_{\rm b}(X,Y)$) the collection of continuous
(and bounded) functionals from $X$ to $Y$.
When $Y=\R$, we simply use $\cC(X)$ ($\cC_{\rm b}(X)$) to represent $\cC(X,\R)$
($\cC_{\rm b}(X,\R)$).

Let $\ell^p$, $p\in[1,\8)$, be the Banach space of all real valued $p$-power summable
bi-infinite sequences with the norm $\|\cdot\|_p$ such that
$$\|u\|_p^p=\sum_{i\in\Z}|u_i|^p,\Hs\mb{for }u=(u_i)_{i\in\Z}\in\ell^p.$$
When $p=\8$, we let $\ell^\8$ be the Banach space of all bounded bi-infinite
sequences with the norm $\|\cdot\|_\8$ satisfying
$$\|u\|_\8=\sup_{i\in\Z}|u_i|,\Hs\mb{for }u=(u_i)_{i\in\Z}\in\ell^\8.$$
For the case when $p=2$, $\ell^2$ is a Hilbert space with the inner product
$$(u,v):=\sum_{i\in\Z}u_iv_i\hs\mb{for all }
u=(u_i)_{i\in\Z},\,v=(v_i)_{i\in\Z}\in\ell^2$$
and the norm simply denoted by $\|\cdot\|$.
These spaces have the following embedding
relationship between $\ell^p$ and $\ell^q$ (see \cite{WJ24}),
\be\label{2.1}1\leqslant p\leqslant q\leqslant\8\hs\Ra\hs \ell^p\subset \ell^q\mb{ and }
\|u\|_q\leqslant\|u\|_p\mb{ for }u\in\ell^p.\ee
For notational convenience in the sequel,
given each $u\in\ell^p$, $p\geqslant1$ and $r>0$,
we use $|u|$, $u^r$ and ${\rm sgn}u$ to denote the sequences as
$$(|u|)_i=|u_i|,\Hs(u^r)_i=u_i^r\Hs\mb{and}\Hs({\rm sgn}u)_i={\rm sgn}u_i,$$
respectively, for each $i\in\Z$,
where ${\rm sgn}:\R\ra\{\pm1,0\}$ is the sign function, i.e.,
${\rm sgn}x=x/|x|$ for $x\neq0$ and ${\rm sgn}0=0$.
According to this, it is easy to see that $|u|^r=(|u_i|^r)_{i\in\Z}$
for $u\in\ell^p$.
When it comes across the calculation $0^0$ in the sequel, we set $0^0$ to be $1$.

Let $u\in\ell^p$, $v\in\ell^q$.
By the embeddings above we have the following inequality that will be used frequently
in the sequel:
\be\label{2.2n}(u,|v|^r)=\sum_{i\in\Z}u_i|v_i|^r\leqslant\sum_{i\in\Z}|u_i|\|v\|_{\8}^r
\leqslant\|u\|_1\|v\|_{q}^r\hs\mb{for }p=1,\,q\geqslant1,\,r\geqslant0.\ee

Define a certain multiplication $\otimes$ for two arbitrary bi-infinite sequences such that
$$(u\otimes v)_i:=u_iv_i\hs\mb{for all }
u=(u_i)_{i\in\Z},\,v=(v_i)_{i\in\Z}.$$
We also define the self-mappings $B$ and $B^*$ on the spaces of bi-infinite sequences
respectively as
$$(Bu)_i:=u_{i+1}-u_i,\hs\mb{and}\hs(B^*u)_i:=u_{i-1}-u_i.$$
Obviously $Au=B^*(Bu)$, $(B^*u,v)=(u,Bv)$.
It is trivial to see that $B$, $B^*$ and $A$ are all bounded linear operators
from $\ell^q$ to $\ell^q$ for all $q\geqslant 1$, i.e.,
\be\label{2.3n}\|Bu\|_q\leqslant2\|u\|_q,\hs\|B^*u\|_q\leqslant2\|u\|_q,
\hs\|Au\|_q\leqslant4\|u\|_q
\hs\mb{and}\hs(A^2u,u)=\|Au\|^2.\ee
Based on the settings above, the equation \eqref{1.1} can be rewritten as
\be\label{2.4n}\di u(t)=-[A(Au+f(u))+\lam u-g(t,u)]\di t+(ku+h(t))\di W(t),\; t>\tau,
\hs u(\tau)=u_\tau\in\ell^2.\ee

\subsection{Nonautonomous random dynamical system}\label{ss2.2}

Consider the space $\cC_0(\R)=\{\w\in\cC(\R):\w(0)=0\}$,
with the compact open topology, the Borel $\sig$-algebra $\cF$ and
the corresponding Wiener measure $\bP$.
For sake of seeking for a solution of \eqref{2.4n},
we need to introduce the OU equation (see \cite{WLYJ21}).
Let $\theta_t\w(\cdot)=\w(\cdot+t)-\w(t)$, $t\in\R$.
Then $(\cC_0(\R),\cF,\bP,\{\theta_t\}_{t\in\R})$ is
an ergodic metric dynamical system.
Based on the system $(\cC_0(\R),\cF,\bP,\{\theta_t\}_{t\in\R})$, we set
the stochastic stationary solution
\be\label{2.5n}z(\theta_t\w)=-\int_{-\8}^0\me^{s}\theta_t\w(s)\di s,\ee
as a pathwise solution of the OU equation
$\di z+z\di t=\di\w(t)$.
Then there is a $\theta_t$-invariant set $\W\subset\cC_0(\R)$ such that
$\bP(\W)=1$ and for every $\w\in\W$,
$z(\theta_t\w)$ is continuous in $t$.
The random variable $z(\theta_t\w)$ satisfies
\be\label{2.6n}
  \lim_{t\ra \pm\infty}\frac{z(\theta_t\w)}{t}
  =\lim_{t\ra \pm\infty}\frac{\int_{0}^{t}z(\theta_s\w)\di s}{t}=0.
\ee
Hence we only consider the set $\W$.
As to the Wiener process $W$ presented in \eqref{1.1},
we define $W$ on the probability space $(\W,\cF,\bP)$ as
$$W(\w)(t):=\w(t)\hs\mb{for an arbitrary }\w\in\W.$$

We now present the definition of NRDS, which is indeed a cocycle of random dynamical systems
(see \cite{WJ24,WZ23}).

\bd\label{de2.1} Let $X$ be a complete metric space and $(\W,\cF,\bP,\{\theta_t\}_{t\in\R})$
be an ergodic metric dynamical system.
A \textbf{nonautonomous random dynamical system} $\vp:\R^+\X\R\X\W\X X\ra X$ is a mapping
such that
\benu\item[(1)] $\vp(\cdot,\tau,\cdot,\cdot):\R^+\X\W\X X\ra X$ is
$(\cB(\R^+)\X\cF\X\cB(X),\cB(X))$-measurable for every $\tau\in\R$;
\item[(2)] $\vp(0,\tau,\w,\cdot)$ is the identity on $X$ for every $(\tau,\w)\in\R\X\W$;
\item[(3)] $\vp(t+s,\tau,\w,\cdot)=\vp(t,\tau+s,\theta_s\w,\vp(s,\tau,\w,\cdot))$
for every $t,s\in\R^+$ and $(\tau,\w)\in\R\X\W$.
\eenu
Moreover, an NRDS $\vp$ is said to be \textbf{continuous}
if $\vp(\cdot,\tau,\w,\cdot):\R^+\X X\ra X$ is continuous for every $(\tau,\w)\in\R\X\W$
for $\bP$-a.s. $\w\in\W$.
\ed

\subsection{Pullback random attractor}

Let $X$ be a Polish space and $2^X$ be the collection of all subsets of $X$.
An $X$-valued time-parametrized random variable $R(\tau,\w)$ is called
a \emph{nonautonomous random variable} in $X$.
A \emph{nonautonomous random set} $D$ in $X$ is a family of nonempty
subsets $D(\tau,\w)\in 2^X$ with two parameters
$(\tau,\w)\in\R\X\W$, which is $(\W,\cF,\bP)$-measurable ($\cF$-measurable for short)
with respect to $\w\in\W$ (see \cite{C02}),
i.e., the mapping $\w\mapsto \di_X(x,D(\tau,\w))$ is $(\cF,\cB(\R))$-measurable for each fixed
$(x,\tau)\in X\X\R$, where $\di_X(\cdot,\cdot)$ is the metric induced by norm of $X$.
If each $D(\tau,\w)$, $(\tau,\w)\in\R\X\W$, is closed (resp. bounded/compact) in $X$,
we say $D$ is a \emph{nonautonomous random closed (resp. bounded/compact) set}.
We have the following result by a simple extension of \cite[Proposition 1.3.6]{C02}.
\bp\label{pro2.2n} Let $V:X\ra\R$ be a continuous function on $X$ and
$R(\tau,\w)$ be a nonautonomous random variable.
If the set $\{x:V(x)\leqslant R(\tau,\w)\}$ is nonempty for each $(\tau,\w)\in\R\X\W$,
then it is a random closed set.\ep

We say $\cD$ is a \emph{universe} in $X$, if $\cD$ is an inclusion-closed collection of nonautonomous
random sets, i.e., when $D\in\cD$ and a nonautonomous random set $D'$ satisfies
$D'(\tau,\w)\subset D(\tau,\w)$ for each $(\tau,\w)\in\R\X\W$, then $D'\in\cD$.

\bd Let $\vp$ be an NRDS on $X$ and $\cD$ a universe in $X$.
A set-mapping $\cA:\R\X\W\ra 2^{X}$ is said to be a
\textbf{pullback $\cD$-attractor} (or \textbf{pullback attractor with respect to} $\cD$)
in $X$ for $\vp$ if
\benu
\item[(1)] $\cA(\tau,\w)$ is compact in $X$ for all $\tau\in\R$ and $\bP$-a.s. $\w\in\W$;
\item[(2)] $\cA$ is invariant under the system $\vp$, i.e.,
      $$\vp(t,\tau,\w,\cA(\tau,\w))=\cA(\tau+t,\theta_t\w),
      \mb{ for all }t\geqslant0\mb{ and }(\tau,\w)\in\R\X\W;$$
\item[(3)] $\cA$ is pullback $\cD$-attracting in $X\cap Y$, that is, for every $D\in\cD$,
      \be\label{2.7n}\lim_{t\ra +\8}{\rm dist}_{X\cap Y}(\vp(t,\tau-t,
      \theta_{-t}\w,D(\tau-t,\theta_{-t}\w)),\cA(\tau,\w))=0;\ee
\item[(4)] there exists a closed set $K\in\cD$ with $\cA\subset K$.
\eenu
The pullback $\cD$-attractor is called an \textbf{pullback random $\cD$-attractor}
in $X$ for $\vp$ if
\benu\item[(5)] $\cA$ is a nonautonomous random set in $X$.
\eenu
\ed

In \eqref{2.7n}, ${\rm dist}_X(\cdot,\cdot)$ is the \emph{Hausdorff semi-distance}
under the norm of an arbitrary normed space $X$, i.e., for two nonempty sets $A,B\subset X$,
$${\rm dist}_X(A,B):=\sup_{a\in A}\inf_{b\in B}\|a-b\|_X.$$

To verify the existence of pullback random attractors,
we also need to recall two related concepts.
Let $\vp$ be an NRDS on $X$ and $\cD$ a universe.
A \emph{pullback $\cD$-absorbing set} for the NRDS $\vp$
is a nonautonomous random closed set $K$ in $X$, such that
for every $D\in\cD$ and all $\tau\in\R$, $\bP$-a.s. $\w\in\W$,
there exists $T=T(\tau,\w,D)>0$ such that for all $t>T$,
$$\vp(t,\tau-t,\theta_{-t}\w,D(\tau-t,\theta_{-t}\w))\subset K(\tau,\w).$$
We say $\vp$ is \emph{pullback $\cD$-asymptotically compact} in $X$,
if for each $\tau\in\R$ and $\bP$-a.s. $\w\in\W$,
the sequence $\vp(t_n,\tau-t_n,\theta_{-t_n}\w,x_n)$
has a convergent subsequence in $Y$ whenever $t_n\ra+\8$ and
$x_n\in D(\tau-t_n,\theta_{-t_n}\w)$ with $D\in\cD$.

For existence of pullback random attractors,
we have the following result (see similar results in \cite{WJ24,WZ23}).

\bt\label{th2.4n}
Suppose that $\vp$ is a continuous NRDS on $X$ over $(\W,\cF,\bP,{\{\theta_t\}}_{t\in\R})$
with a universe $\cD$ in $X$. Assume that
\benu
\item[(1)] $\vp$ has a pullback $\cD$-absorbing set $K$ in $X$;
\item[(2)] $\vp$ is pullback $\cD$-asymptotically compact in $X$.
\eenu
Then $\vp$ has a unique pullback random $\cD$-attractor $\cA$ in $X$.
Moreover, if $K\in\cD$, then $\cA\in\cD$.
\et

\subsection{Invariant sample measures}

We first give the definition of invariant sample measures (with two variables) for NRDS's;
see \cite{WZ23,CY23}.
\bd Let $\vp$ be an NRDS on a Banach space $X$.
A family of Borel probability measures $\{\mu_{\tau,\w}\}_{(\tau,\w)\in\R\X\W}$ on $X$ is called
the \textbf{invariant sample measures} for $\vp$, if for each $t\geqslant0$, $\tau\in\R$,
$\bP$-a.s. $\w\in\W$ and all $U\in\cB(X)$,
$$
\mu_{\tau+t,\theta_t\w}(U)=\mu_{\tau,\w}((\vp(t,\tau,\w,\cdot))^{-1}U).
$$
\ed

In order to construct the invariant sample measures, we need to introduce
the \emph{generalized Banach limit}, which is any linear functional,
denoted by $\disp\LIM_{t\ra-\8}$ (we only need to consider the case $t\ra-\8$,
see more in \cite{WZC20,WZ23}), defined for an arbitrary bounded real-valued function
on $(-\8,a]$ for some $a\in\R$ and satisfying
\benu
\item[(1)] $\disp\LIM_{t\ra-\8}\zeta(t)\geqslant0$ for nonnegative functions $\zeta$
on $(-\8,a]$;
\item[(2)] $\disp\LIM_{t\ra-\8}\zeta(t)=\lim_{t\ra-\8}\zeta(t)$ if the latter limit exists.
\eenu

For the existence of invariant sample measures for NRDS's,
one can refer to Theorem 2.6 in \cite{WZ23} as follows.

\bp\label{pro2.6n}
Let $(X,d)$ be a complete metric space, $\vp$ be an NRDS over
$(\W,\cF,\bP,\{\theta_t\}_{t\in\R})$ on the state space $(X,d)$ and $\cD$ be
a universe for $\vp$. Suppose that
\benu\item[(1)] $\vp$ has a nonautonomous random $\cD$-attractor
$\{\cA(\tau,\w)\}_{(\tau,\w)\in\R\X\W}$ in $X$;
\item[(2)] for each $t,\,\tau\in\R$ with $t\geqslant\tau$ and $\bP$-a.s. $\w\in\W$,
the $X$-valued mapping
$$(\tau,u)\mapsto\vp(t-\tau,\tau,\theta_{\tau}\w,u)$$
is continuous and bounded for $\tau\in(-\8,t]$ and each fixed $u\in X$.\eenu
Then given a generalized Banach limit $\disp\LIM_{\tau\ra-\8}$ and a
mapping $\xi:\R\X\W\ra X$ with $\{\xi(\tau,\w)\}_{(\tau,\w)\in\R\X\W}\in\cD$
such that $\xi(\tau,\theta_\tau\w)$ is continuous in $\tau$,
there exists a unique family of Borel probability measures
$\{\mu_{\tau,\w}\}_{(\tau,\w)\in\R\X\W}$ in $X$ such that the support of
$\mu_{\tau,\w}$ is contained in $\cA(\tau,\w)$ and for all $\Upsilon\in\cC(X)$,
$$\ba{rl}&\disp\LIM_{\tau\ra-\8}\frac1{t-\tau}\int_{\tau}^t
 \Upsilon(\vp(t-s,s,\theta_s\w,\xi(s,\theta_s\w)))\di s
 =\int_{\cA(t,\theta_t\w)}\Upsilon(u)\di\mu_{t,\theta_t\w}(u)\\
=&
 \disp\int_X\Upsilon(u)\di\mu_{t,\theta_t\w}(u)
 =\LIM_{\tau\ra-\8}\frac1{t-\tau}\int_{\tau}^t\int_X
 \Upsilon(\vp(t-s,s,\theta_s\w,u))\di\mu_{s,\theta_s\w}(u)\di s.\ea
$$
Moreover, the measure $\{\mu_{\tau,\w}\}_{(\tau,\w)\in\R\X\W}$ is invariant in the sense that
$$\int_{\cA(\tau+t,\theta_t\w)}\Upsilon(u)\di\mu_{\tau+t,\theta_t\w}(u)
=\int_{\cA(\tau,\w)}\Upsilon(\vp(t,\tau,\w,u))
\di\mu_{\tau,\w}(u),\hs\mb{for all }t\geqslant0.$$
\ep
\subsection{Sample statistical solutions}\label{ss2.5}

We develop in this subsection the sample statistical solutions with respect to
the invariant sample measures given in Proposition \ref{pro2.6n}.

In Section \ref{s1}, we have given the class of test functions in Definition \ref{de1.1}
and the definition of the stochastic Liouville type equation.
Here we first present some remarks on these concepts.

\br The class $\cT$ in Definition \ref{de1.1}
is not empty for \eqref{1.5} with suitable solutions
(see the specific case in \eqref{3.6z} below for details).
Actually Definition \ref{de1.1} is an extension of the kinds of test functions defined
as the following example (see \cite[Chapter V]{FMRT01}).
For example, we can take the cylindrical functions $\Psi:H\ra\R$ in $\cT$ of the following form
$$\Psi(u):=\psi((u,e_1),\cdots,(u,e_m)),$$
where $\psi$ is a $C_c^2$ scalar valued function defined on $\R^m$, $m\in\N$
and $e_1$, $\cdots$, $e_m$ belong to $V$ given in Section \ref{s1}
(see \cite{CY23,FMRT01} for more details).
Then we know for all $v,w\in H$,
$$\(\Psi'(u),v\)=\sum_{i=1}^m\pa_i\psi\((u,e_1),\cdots,(u,e_m)\)(e_i,v),$$
$$\Psi''(u)(v,w)=\sum_{i,j=1}^m\pa_{ij}\psi\((u,e_1),\cdots,(u,e_m)\)(e_i,v)(e_j,w),$$
where $\pa_i\psi$ is the derivative of $\psi$ with respect to the $i$-th variable,
and $\pa_{ij}\psi$ is the second-order derivative of $\psi$ with respect to
the $i$-th and $j$-th variables.
\er
\br Let $\{\mu_{\tau,\w}\}_{(\tau,\w)\in\R\X\W}$ be a family of invariant sample measures given
by Proposition \ref{pro2.6n}.
We say the invariant sample measures $\{\mu_{\tau,\w}\}_{(\tau,\w)\in\R\X\W}$ satisfy
the \textbf{stochastic Liouville type theorem}, if the stochastic Liouville type equation \eqref{1.7} holds
for all $t\geqslant s$.

Dividing \eqref{1.7} by $t-s$ and letting $t$ tend to $s$,
we can obtain by the invariance of $\mu_{\tau,\w}$
the differential form of stochastic Liouville type equation as follows,
\begin{align*}\di\int_{H}\Psi(u)\mu_{s,\theta_{s-\tau}\w}(\di u)
=&\int_H\langle\tilde{F}(s,u),\Psi'(u)\rangle
\mu_{s,\theta_{s-\tau}\w}(\di u)\di s\\
&+\int_H(\tilde{G}(s,u),\Psi'(u))\mu_{s,\theta_{s-\tau}\w}(\di u)
\di\tilde{W}(s)\\
&+\frac12\int_H\Psi''(u)(\tilde{G}(s,u),\tilde{G}(s,u))
\mu_{s,\theta_{s-\tau}\w}(\di u)\di s.
\end{align*}
\er
Based on the extended stochastic Liouville type equation \eqref{1.7},
we can further develop the statistical solutions for deterministic system into the
stochastic situation, where, we call them ``sample statistical solutions".
We follow the class $\cT$ defined in Definition \ref{de1.1} to
give the definition of sample statistical solutions.

\bd\label{de2.9n} A family $\{\mu_{\tau,\w}\}_{(\tau,\w)\in\R\X\W}$ of
Borel probability measures in $H$ is called a \textbf{sample statistical solution}
(in $H$) of \eqref{1.5}
if the following conditions are satisfied:
\benu\item[(1)] The function $s\mapsto\int_H\Psi(u)\mu_{s,\theta_{s-\tau}\w}(\di u)$
is continuous for every $\Psi\in\cC_b(H)$;
\item[(2)] For almost all $s\in\R$, the functions
$$u\mapsto\langle\tilde{F}(s,u),\phi\rangle,\hs u\mapsto(\tilde{G}(s,u),\phi)
\hs\mb{and}\hs u\mapsto\Phi(\tilde{G}(s,u),\tilde{G}(s,u))$$
are $\mu_{s,\theta_{s-\tau}\w}$-integrable for every $\phi\in V$ and bilinear mapping
$\Phi\in L(H\X H,\R)$;
moreover, the mappings
$$s\mapsto\int_H\langle\tilde{F}(s,u),\phi\rangle\mu_{s,\theta_{s-\tau}\w}(\di u)
\hs\mb{and}\hs
s\mapsto\int_H\Phi(\tilde{G}(s,u),\tilde{G}(s,u))\mu_{s,\theta_{s-\tau}\w}(\di u)$$
belongs to $L_{\rm loc}^1(\R)$, and the mapping
$$s\mapsto\int_H\(\tilde{G}(s,u),\phi\)\mu_{s,\theta_{s-\tau}\w}(\di u)$$
belongs to $L^2_{\rm loc}(\R)$;
\item[(3)] For each $\Psi\in\cT$, the stochastic Liouville type equation
\eqref{1.7} holds.
\eenu
\ed

\section{The case when $k>0$}

In this section we consider the problem \eqref{2.4n} for the case when $k>0$.
For relative calculations, we first require some hypotheses
for the mappings $f$, $g$ and $h$ in \eqref{2.4n} as follows.

\noindent\textbf{(F)} The mapping $f$ is locally Lipschitz, and there are $q>2$, $\epsilon>0$
and $\alpha_{k}=(\alpha_{ki})_{i\in\Z}$, $k=1,2$, with ${q-1-\epsilon}>0$, such that
\be|f_i(u_i)-f_i(v_i)|\leqslant\gamma(R)|u_i-v_i|,\hs\mb{whenever }|u_i|,\,|v_i|\leqslant R\label{3.1z}\ee
$$\mb{and}\hs|f_i(u_i)|\leqslant\alpha_{1i}|u_i|^{q-1-\epsilon}+\alpha_{2i}\,
\mb{ with }\alpha_{1},\,\alpha_{2}\in\ell^1,$$
where $\gam:\R^+\ra\R^+$ is such that
$\gam(\cI)$ is bounded for each bounded interval $\cI\subset\R^+$.

\noindent\textbf{(G1)} The mapping $(t,u_i)\mapsto g_i(t,u_i)$ is continuous
in both variables and differentiable in $u_i$,
and there are $\beta>0$ and $\psi_k(t)=(\psi_{ki}(t))_{i\in\Z}$, $k=1,
2,3,4$, with $q>2$ given in \textbf{(F)} such that
\begin{align}\label{3.2z}g_i( t,u_i)u_i
&\leqslant-\beta|u_i|^q+\psi_{1i}(t),
\mb{ with }\psi_1\in L_{\rm loc}^1(\R;\ell^1),\\
\label{3.3z}\frac{\pa g_i}{\pa u_i}(t,u_i)
&\leqslant\psi_{2i}(t),\mb{ with }\psi_2\in L_{\rm loc}^\infty(\R;\ell^\infty)\hs\mb{and}\\
|g_i(t,u_i)|&\leqslant\psi_{3i}(t)|u_i|^{q-1}+\psi_{4i}(t),
\mb{ with }\psi_3\in L_{\rm loc}^{\infty}(\R,\ell^{\infty})
\mb{ and }\psi_4\in L_{\rm loc}^{2}(\R,\ell^{2}).
\label{3.4z}\end{align}

\noindent\textbf{(H1)} The mapping $h:\R\ra\ell^2\cap\ell^q$ ($q>2$ given in \textbf{(F)})
is such that $h:\R\ra\ell^2$ is derivable,
\be\label{3.5z}h\in L^2_{\rm loc}(\R,\ell^2)\cap L^q_{\rm loc}(\R,\ell^q)\hs\mb{and}
\hs h'\in L^2_{\rm loc}(\R,\ell^2).\ee

\subsection{Unique existence of solutions}\label{ss3.1}

A process $\{u(t,\tau,\w,u_{\tau})\}_{t\geqslant\tau,\w\in\W}$
is called a \emph{local solution} of problem \eqref{2.4n}, provided there exists
$T>\tau$ such that for $(\tau,u_\tau)\in\R\X\ell^2$ and $\bP$-a.s. $\w\in\W$,
$$
  u(\cdot,\tau,\w,u_{\tau})\in\cC^1([\tau,T);\ell^2)\cap L^q(\tau,T;\ell^q)
$$
and $u(\cdot,\tau,\w,u_\tau)$ satisfies \eqref{2.4n}.
If additionally for $(\tau,u_\tau)\in\R\X\ell^2$ and $\bP$-a.s. $\w\in\W$,
\be\label{3.6z}
  u(\cdot,\tau,\w,u_{\tau})\in\cC^1([\tau,\8);\ell^2)\cap L^q_{\rm loc}(\tau,\8;\ell^q),
\ee
we say the solution $u(t,\tau,\w,u_{\tau})$ is \emph{global}.
Here we use $\cC^1$ to mean differentiability and continuity of the derivative.

Similar to the methods applied in \cite{WLYJ21,WZ23},
we transform the stochastic equation \eqref{1.1} into a random one,
so that the solution \eqref{2.4n} can generate a nonautonomous random dynamical system.
Let $u(t,\tau,\w,u_\tau)$ be a solution of \eqref{2.4n} and set
\be\label{3.7z}v(t)=\me^{-kz(\theta_t\w)}(u(t)+k^{-1}h(t)),\Hs
v(\tau)=v_\tau:=\me^{-kz(\theta_\tau\w)}(u_\tau+k^{-1}h(\tau)),\ee
where $z(\theta_t\w)$ is given in \eqref{2.5n}.
We deduce from \eqref{2.4n} and \eqref{3.7z} that
$v(t)$ satisfies the following random system
\be\label{3.8z}
\frac{\di v(t)}{\di t}=F(t,\w,v),\hs t>\tau,\Hs v(\tau)=v_\tau,\ee
where
\begin{align*}F(t,\w,v):=&\(-A^2v+\(kz(\theta_t\w)-\lam\)v\)\\
&+\me^{-kz(\theta_t\w)}\(-Af(\me^{kz(\theta_t\w)}v-k^{-1}h(t))
+g(t,\me^{kz(\theta_t\w)}v-k^{-1}h(t))\)\\
&+k^{-1}\me^{-kz(\theta_t\w)}\(A^2h(t)+\lam h(t)+h'(t)\)\\
:=&F_1(t,\w,v)+F_2(t,\w,v)+F_3(t,\w).
\end{align*}

The unique existence theorem below of solutions is summarized in \cite{WJ24}.
\bt\label{th3.1} Let $X$ be a real Hilbert space with inner product $(\cdot,\cdot)$
and norm $\|\cdot\|$,
$\cI=[0,a]\subset\R$, $E=\{x\in X:\|x-x_0\|\leqslant r\}$ for some $r>0$,
$f:\cI\X E\ra X$ continuous and $|f(t,x)|\leqslant m$ on $\cI\X E$.
Let $f$ satisfy the condition
$$(f(t,x)-f(t,y),x-y)\leqslant C\|x-y\|^2\hs\mb{for }t\in(0,a],\;x,y\in E$$
for some positive constant $C$ independent of $t$, $x$ and $y$.
Let $b<\min\{a,r/m\}$.
Then the differential equation $x'=f(t,x)$, $x(0)=x_0$ has a unique solution $x(t;x_0)$ on $[0,b]$
such that $x\in\cC^1([0,b];X)$.
\et

In the following argument, for notational convenience, we always denote $c$ as
an arbitrary positive constant only depending on $\lam$, $k$, $q$, $\epsilon$ and $\beta$
given in the hypotheses \textbf{(F)},
\textbf{(G1)} and \textbf{(H1)}
and $c$ may be different from line to line and even in the same line.

Based on Theorem \ref{th3.1}, we can obtain the solution of the system \eqref{3.8z}
in the following theorem.
\bt\label{th3.2} Let the hypotheses \textbf{(F)}, \textbf{(G1)} and \textbf{(H1)} hold.
Then for each $(\tau,v_\tau)\in\R\X\ell^2$ and $\bP$-a.s. $\w\in\W$,
the system \eqref{3.8z} has a unique global solution
$v(\cdot,\tau,\w,v_\tau)\in\cC^1([\tau,\8);\ell^2)\cap L^q_{\rm loc}(0,\8;\ell^q)$
such that $v(\tau,\tau,\w,v_\tau)=v_\tau$.
\et
\bo We first use Theorem \ref{th3.1} to discuss the local unique existence
of solutions of $\di v=F(t,\w,v)\di t$.
Let $T>\tau$ and $t\in[\tau,T]$ and pick $u,v\in\ell^2$ such that $\|u\|$, $\|v\|<R$.
By \eqref{2.3n}, we know
\be\label{3.9z}
\(F_1(t,\w,u)-F_1(t,\w,v),u-v\)=-\|A(u-v)\|^2+(kz(\theta_t\w)-\lam)\|u-v\|^2.
\ee
By \eqref{3.1z}, we have
\begin{align}&\(A(f(\me^{kz(\theta_t\w)}u-k^{-1}h(t))
-f(\me^{kz(\theta_t\w)}v-k^{-1}h(t))),u-v\)\notag\\
=&\(f(\me^{kz(\theta_t\w)}u-k^{-1}h(t))-f(\me^{kz(\theta_t\w)}v-k^{-1}h(t)),A(u-v)\)\notag\\
 \leqslant&\|f(\me^{kz(\theta_t\w)}u-k^{-1}h(t))
 -f(\me^{kz(\theta_t\w)}v-k^{-1}h(t))\|\|A(u-v)\|\notag\\
\leqslant&4\gamma(R)\me^{kz(\theta_t\w)}\|u-v\|^2.
 	\label{3.10z}\end{align}
For the nonlinear term $g$, by \eqref{3.3z} we similarly have
\begin{align}&\(g(t,\me^{kz(\theta_t\w)}u-k^{-1}h(t))
-g(t,\me^{kz(\theta_t\w)}v-k^{-1}h(t)),u-v\)\notag\\
\leqslant&\me^{kz(\theta_t\w)}\sum_{i\in\Z}\psi_{2i}(t)|u_i-v_i|^2
\leqslant\me^{kz(\theta_t\w)}\|\psi_2(t)\|_\8\|u-v\|^2.
\label{3.11z}\end{align}

We can then deduce by \eqref{2.3n}, \eqref{3.9z}, \eqref{3.10z} and \eqref{3.11z} that
\begin{align}&\(F(t,\w,u)-F(t,\w,v),u-v\)\notag\\
\leqslant&\(\esssup_{t\in[\tau,T]}\|\psi_2(t)\|_\8+4\gamma(R)
+\max_{t\in[\tau,T]}|kz(\theta_t\w)-\lam|\)\|u-v\|^2\notag\\
:=&\cP_1(\tau,T,\w,R)\|u-v\|^2.\label{3.12z}
\end{align}
Then the local unique existence of solutions for
the problem \eqref{3.8z} is derived from Theorem \ref{th3.1},
i.e., for each $(\tau,v_\tau)\in\R\X\ell^2$
and $\bP$-a.s. $\w\in\W$, \eqref{3.8z} has a solution
$v(\cdot,\tau,\w,v_\tau)\in\cC^1([\tau,T);\ell^2)$
for some $T\in(\tau,+\8)$.

In order to show the global existence of the solution, we need another estimate in advance.
First, we observe that
\be\label{3.13z}(F_1(t,\w,v),v)=-\|Av\|^2+(kz(\theta_t\w)-\lam)\|v\|^2\hs\mb{and}\ee
\begin{align}(F_2(t,\w,v),v)=&\(F_2(t,\w,v),v-k^{-1}\me^{-kz(\theta_t\w)}h(t)\)
+\(F_2(t,\w,v),k^{-1}\me^{-kz(\theta_t\w)}h(t)\)\notag\\
:=&I_1+I_2.\end{align}
By \eqref{2.2n}, the hypothesis \textbf{(F1)}, Young's and H\"older's inequality, we see that
\begin{align}\me^{2kz(\theta_t\w)}I_1\leqslant&
\left|\(f(\me^{kz(\theta_t\w)}v-k^{-1}h(t)),A\(\me^{kz(\theta_t\w)}v-k^{-1}h(t)\)\)\right|\notag\\
&+\(g(t,\me^{kz(\theta_t\w)}v-k^{-1}h(t)),\me^{kz(\theta_t\w)}v-k^{-1}h(t)\)\notag\\
\leqslant&\sum_{i\in\Z}\(\al_{1i}|\me^{kz(\theta_t\w)}v_i-k^{-1}h_i(t)|^{q-1-\epsilon}+\al_{2i}\)
\left|\(A(\me^{kz(\theta_t\w)}v-k^{-1}h(t))\)_i\right|\notag\\
&-\beta\|\me^{kz(\theta_t\w)}v-k^{-1}h(t)\|_q^q+\|\psi_1(t)\|_1\notag\\
\leqslant&4\(\|\al_1\|_1\|\me^{kz(\theta_t\w)}v-k^{-1}h(t)\|_q^{q-\epsilon}
+\|\al_2\|_1\|\me^{kz(\theta_t\w)}v-k^{-1}h(t)\|_q\)\notag\\
&-\beta\|\me^{kz(\theta_t\w)}v-k^{-1}h(t)\|_q^q+\|\psi_1(t)\|_1\notag\\
\leqslant&\frac{4\epsilon}{q}
\(\frac{32(q-\epsilon)}{\beta q}\)^{\frac{q-\epsilon}{\epsilon}}\|\al_1\|_1^{\frac{q}{\epsilon}}
+\frac{4}{p}\(\frac{32}{\beta q}\)^{\frac1{q-1}}\|\al_2\|_1^{\frac{q}{q-1}}\notag\\
&-\frac{3\beta}{4}\|\me^{kz(\theta_t\w)}v-k^{-1}h(t)\|_q^q+\|\psi_1(t)\|_1\notag\\
\leqslant&c\(\|\al_1\|_1^{\frac{q}{\epsilon}}+\|\al_2\|_1^{\frac{q}{q-1}}\)+\|\psi_1(t)\|_1
-\frac{3\beta}{4}\|\me^{kz(\theta_t\w)}v-k^{-1}h(t)\|_q^q
\end{align}
and similarly
\begin{align}\me^{2kz(\theta_t\w)}|I_2|\leqslant&
\left|\(f(\me^{kz(\theta_t\w)}v-k^{-1}h(t)),k^{-1}Ah(t)\)\right|
+\left|\(g(t,\me^{kz(\theta_t\w)}v-k^{-1}h(t)),k^{-1}h(t)\)\right|\notag\\
\leqslant&4k^{-1}\(\|\al_1\|_1\|\me^{kz(\theta_t\w)}v-k^{-1}h(t)\|_q^{q-1-\epsilon}
+\|\al_2\|_1\)\|h(t)\|_q\notag\\
&+k^{-1}\sum_{i\in\Z}|\psi_{3i}(t)||\me^{kz(\theta_t\w)}v_i-k^{-1}h_i(t)|^{q-1}|h_i(t)|
+k^{-1}\|\psi_4(t)\|\|h(t)\|\notag\\
\leqslant&c\(\|\al_1\|_1^{\frac{q}{\epsilon}}+\|\al_2\|_1^{\frac{q}{q-1}}
+\|\psi_4(t)\|^2+\|h(t)\|^2+\(1+\|\psi_3(t)\|^q_\8\)\|h(t)\|_q^q\)\notag\\
&+\frac{\beta}{4}\|\me^{kz(\theta_t\w)}v-k^{-1}h(t)\|_q^q.
\end{align}
Moreover, we also have that
\be\label{3.17z}|(F_3(t,\w),v)|\leqslant
c\me^{-2kz(\theta_t\w)}\(\|h(t)\|^2+\|h'(t)\|^2\)+\frac{\lam}{2}\|v\|^2.
\ee

Now taking the inner product of \eqref{3.8z} and $v$,
we obtain by the estimates from \eqref{3.13z} to \eqref{3.17z} that
\be\frac{\di\|v(t)\|^2}{\di t}+(\lam-2kz(\theta_t\w))\|v\|^2
+\beta\me^{-2kz(\theta_t\w)}\|\me^{kz(\theta_t\w)}v-k^{-1}h(t)\|_q^q
\leqslant c\me^{-2kz(\theta_t\w)}\cQ_1(t),\label{3.18z}
\ee
where
\begin{align}\notag\cQ_1(t):=&\|\al_1\|_1^{\frac{q}{\epsilon}}+\|\al_2\|_1^{\frac{q}{q-1}}
+\|\psi_1(t)\|_1+\|\psi_4(t)\|^2\\
&+\(1+\|\psi_3(t)\|^q_\8\)\|h(t)\|_q^q+\|h(t)\|^2+\|h'(t)\|^2.\label{3.19z}\end{align}
Applying Gronwall's lemma to \eqref{3.18z} over the time interval $[\sig,s]\subset[\tau,T)$,
we obtain
\begin{align}&\|v(s)\|^2+\beta\int_\sig^s\me^{\lam(\vsig-s)+2k\int_\vsig^sz(\theta_r\w)\di r
-2kz(\theta_\vsig\w)}\|\me^{kz(\theta_\vsig\w)}v(\vsig)-k^{-1}h(\vsig)\|_q^q\di\vsig\notag\\
\leqslant&\me^{\lam(\sig-s)+2k\int_\sig^sz(\theta_r\w)\di r}\|v(\sig)\|^2
+c\int_{\sig}^s\me^{\lam(\vsig-s)+2k\int_\vsig^sz(\theta_r\w)\di r
-2kz(\theta_\vsig\w)}\cQ_1(\vsig)\di\vsig.
\label{3.20z}\end{align}
Replacing $\sig$ and $s$ by $\tau$ and $t$ respectively for $t\in[\tau,T)$
and by the hypotheses \textbf{(F)}, \textbf{(G1)} and \textbf{(H1)}, we have
\begin{align}&\|v(t)\|^2+\beta\int_\tau^t\me^{\lam(\vsig-t)+2k\int_\vsig^tz(\theta_r\w)\di r
-2kz(\theta_\vsig\w)}\|\me^{kz(\theta_\vsig\w)}v(\vsig)-k^{-1}h(\vsig)\|_q^q\di\vsig\notag\\
\leqslant&\me^{2k\int_\tau^T|z(\theta_r\w)|\di r}\|v_\tau\|^2
+c\int_{\tau}^T\me^{2k\int_\vsig^T|z(\theta_r\w)|\di r-2kz(\theta_\vsig\w)}\cQ_1(\vsig)\di\vsig\notag\\
:=&\cP_2:=\cP_2(\tau,T,\w,\|v_\tau\|)<\8.\label{3.21z}
\end{align}

Next we extend the local solution to a global one.
Indeed, we only need to show that for each $T\in(\tau,\8)$,
if there is a unique solution $v(t)$ on bounded interval $[\tau,T)$,
then for each increasing sequence $t_n\in[\tau,T)$ with $t_n\ra T^-$,
$v(t_n)$ converges in $\ell^2$,
which implies that each of such sequences $v(t_n)$ possesses the same limit $v^*$
(similar argument can be seen in \cite{WJ24}).
As a result, it suffices to show that $v(t_n)$ is a Cauchy sequence in $\ell^2$.
Consider $\|v(t+\De t)-v(t)\|$ with $\tau<t<t+\De t<T$ and $\De t$ fixed.
By \eqref{3.12z}, we know that
$$\frac{\di}{\di t}\|v(t+\De t)-v(t)\|^2\leqslant
2\cP_1(\tau,T,\w,\sqrt{\cP_2})\|v(t+\De t)-v(t)\|^2.
$$
Then it is obvious that
$$\|v(t+\De t)-v(t)\|\leqslant\|v(\tau+\De t)-v(\tau)\|
\me^{2\cP_1(\tau,T,\sqrt{\cP_2},\w)(T-\tau)}.$$
Pick $t=t_m$ and $\De t=t_n-t_m$ with $n>m$ and we get
\be\label{3.22z}\|v(t_n)-v(t_m)\|\leqslant\|v(\tau+t_n-t_m)-v_\tau\|
\me^{2\cP_1(\tau,T,\sqrt{\cP_2},\w)(T-\tau)}.\ee
Since $t_n$ is increasing, and $t_n\ra T^-$,
one can see that $v(t_n)$ is a Cauchy sequence in $\ell^2$
by applying the continuity of $v(t)$ over $[\tau,T)$ to \eqref{3.22z}.

At last, we prove that $v(\cdot,\tau,\w,v_\tau)\in L^q_{\rm loc}(0,\8;\ell^q)$.
Actually, given $T\in(\tau,\8)$, we take
$$\cP_3(\tau,T,\w):=\exp\(\min_{s\in[\tau,T]}\(
2k\int_s^Tz(\theta_r\w)\di r+k(q-2)z(\theta_s\w)\)+\lam(\tau-T)\).$$
Then following \eqref{3.20z} and \eqref{3.21z}, we know that
\begin{align*}\int_{\tau}^T\|v(s)\|_q^q\di s\leqslant
&2^q\(\int_{\tau}^T\|v(s)-k^{-1}\me^{-kz(\theta_s\w)}h(s)\|_q^q\di s
+k^{-q}\int_\tau^T\|\me^{-kz(\theta_s\w)}h(s)\|_q^q\di s\)\\
\leqslant&\frac{2^q}{\cP_3(\tau,T,\w)}\int_{\tau}^{T}\me^{\lam(s-T)
+2k\int_s^Tz(\theta_r\w)\di r-2kz(\theta_s\w)}\|\me^{kz(\theta_s\w)}v(s)-k^{-1}h(s)\|_q^q\di s\\
&+\(\frac{2}{k}\)^q\int_\tau^T\|\me^{-kz(\theta_s\w)}h(s)\|_q^q\di s<\8,\end{align*}
by \textbf{(H1)}, \eqref{3.21z} and the continuity of $z(\theta_t\w)$ in $t$.
The proof is complete.
\eo

We have obtained the unique existence of the global solution $v(t,\tau,\w,v_\tau)$
for \eqref{3.8z} from Theorem \ref{th3.2}.
Set
\be\label{3.23z}u(t,\tau,\w,u_\tau)=\me^{kz(\theta_t\w)}v(t,\tau,\w,
\me^{-kz(\theta_\tau\w)}(u_\tau+k^{-1}h(\tau)))-k^{-1}h(t).\ee
As a result, $u(t,\tau,\w,u_\tau)$ is exactly the unique global solution of \eqref{2.4n}
for each $(\tau,u_\tau)\in\R\X\ell^2$ and $\bP$-a.s. $\w\in\W$.

\subsection{The induced nonautonomous random dynamical system}\label{ss3.2}

According to Theorem \ref{th3.2},
we define a mapping $\vp:\R^+\X\R\X\W\X\ell^2\ra\ell^2$,
such that for every $(t,\tau,\w,u_\tau)\in\R^+\X\R\X\W\X\ell^2$,
\be\label{3.24z}\vp(t,\tau,\w,u_\tau)=u(t+\tau,\tau,\theta_{-\tau}\w,u_\tau)
=\me^{kz(\theta_t\w)}v(t+\tau,\tau,\theta_{-\tau}\w,v_\tau)-k^{-1}h(t+\tau)\ee
with $v_\tau=\me^{-kz(\theta_\tau\w)}(u_\tau+k^{-1}h(\tau))$.
In order to guarantee that the mapping $\vp$ defined as above
is a continuous NRDS on $\ell^2$ over the system $\(\W,\cF,\bP,\{\theta_t\}_{t\in\R}\)$,
we still need to show the continuity and measurability as in the following theorem.
\bt\label{th3.3} Let the hypotheses \textbf{(F)}, \textbf{(G1)} and \textbf{(H1)} hold
and $u(t,\tau,\w,u_\tau)$ be a solution of \eqref{2.4n}.
Then the following results hold:
\benu\item[(1)]for each $\tau\in\R$ and $\bP$-a.s. $\w\in\W$,
the mapping $(t,u_\tau)\mapsto u(t,\tau,\w,u_\tau)$ is continuous;
\item[(2)]for each $\tau\in\R$,
the mapping $(t,\w,u_\tau)\mapsto u(t,\tau,\w,u_\tau)$ is
$(\cB([\tau,\8))\X\cF\X\cB(\ell^2),\cB(\ell^2))$-measurable.\eenu
\et
\bo For the conclusion, it is sufficient to consider the continuity of $v(t,\tau,\w,v_\tau)$
with respect to $(t,v_\tau)$ according to the relationship \eqref{3.7z}.
Let $v(t)$ and $w(t)$ be solutions of \eqref{3.8z} with respect to the initial data
$v_\tau$ and $w_\tau$, respectively.
Considering the inner product of $\frac{\di(v-w)}{\di t}$ and $v-w$
and by \eqref{3.8z} and \eqref{3.12z}, we obtain
$$\frac12\frac{\di\|v(t)-w(t)\|^2}{\di t}
\leqslant\cP_1\(\tau,T,\w\)\|v-w\|^2.$$
And then for each $T>\tau$, when $t\in[\tau,T]$,
$$\|v(t)-w(t)\|^2\leqslant\me^{2(T-\tau)\cP_1\(\tau,T,\w\)}\|v_\tau-w_\tau\|^2.
$$
This indicates the continuity of the mapping $v_\tau\mapsto v(t,\tau,\w,v_\tau)$ is uniform
for $t\in[\tau,T]$,
which ensures the continuity of $(t,v_\tau)\mapsto v(t,\tau,\w,v_\tau)$
by combining the continuity of $t\mapsto v(t,\tau,\w,v_\tau)$ (given in Theorem \ref{th3.2}).
The details are referred to the proof of \cite[Lemma 4.5]{WZ23}.

The measurability in (2) is a trivial consequence
by the proof procedure of local existence for the solutions of \eqref{3.8z}
in \cite[Theorem 3.2]{D77} and Carath\'eodory Theorem
(one can refer to the proof of \cite[Theorem 2.5]{WJ24} if necessary).
The proof is finished.
\eo

\subsection{Existence of pullback random attractors}\label{s3.2}
We first present the framework for the pullback random attractors.
For the following study, we introduce two arbitrary constants $\lam_0,\,\lam_1\in\R$ such that
$$0<\lam_1<\lam_0<\lam.$$
In the sequel, we always choose $\lam_0$ and $\lam_1$ to be fixed as above.

We now give the further hypotheses on $g$ and $h$ for the existence of pullback random attractors.

\noindent\textbf{(G2)} The functions $\psi_i$ ($i=1,\,3,\,4$) given in \textbf{(G1)} satisfy
\be\label{3.25z}
\int_{-\8}^0\me^{\lam_1 s}\(\|\psi_1(s)\|_1+\|\psi_4(s)\|^2\)\di s<+\8
\hs\mb{and}\hs\psi_3\in L^{\8}((-\8,0],\ell^{\8}).\ee

\noindent\textbf{(H2)} The function $h(t)$ satisfies that
\be\int_{-\8}^{0}\me^{\lam_1 s}\(\|h(s)\|^2+\|h(s)\|_q^q+\|h'(s)\|^2\)\di s<+\8.
\label{3.26z}\ee

There are some simple conclusions deduced from \eqref{3.25z} and \eqref{3.26z} as follows:
\benu\item[(1)] The statements \eqref{3.25z} and \eqref{3.26z} also hold true
if each number $0$ appearing therein is replaced by $t\in\R$,
according to \eqref{3.2z}, \eqref{3.4z} and \eqref{3.5z};
\item[(2)] Each term of the integrands in \eqref{3.25z} and \eqref{3.26z} tends to $0$,
as $s\ra-\8$, e.g.,
\be\label{3.27z}\lim_{s\ra-\8}\me^{\lam_1s}\|h(s)\|^2=0.\ee
\eenu

Let $D=\{D(\tau,\w):\tau\in\R,\w\in\W\}$ be a nonautonomous random set in
$\ell^2$ such that for every $\tau\in\R$ and $\w\in\W$,
\be\label{3.28z}\lim_{s\ra -\8}\me^{\lam_0s+2k\int_{s}^0z(\theta_{r}\w)\di r
-2kz{(\theta_{s}\w)}}\|D(\tau+s,\theta_{s}\w)\|^2=0,\ee
where $\|D\|=\sup_{u\in D}\|u\|$.
In the following discussion of this section,
we always let $\cD$ be the universe of all the nonautonomous random sets $D$
satisfying \eqref{3.28z}.

In this subsection we consider the existence of pullback random attractors
in $\ell^2$.
First we give the existence of pullback random absorbing sets in the following lemma.
\bl\label{le3.4} Let the hypotheses \textbf{(F)}, \textbf{(G1)}, \textbf{(G2)},
\textbf{(H1)} and \textbf{(H2)} hold and
$\me^{kz(\theta_{-t}\w)}v_{\tau-t}-k^{-1}h(\tau-t)\in D(\tau-t,\theta_{-t}\w)$ with $D\in\cD$.
Then for $\bP$-a.s. $\w\in\W$, there exists a certain $T_0=T_0(\tau,\w,D)>0$ such that
when $t>T_0$, it holds that
\begin{align}&\|v(\tau,\tau\!-\!t,\theta_{-\tau}\w,v_{\tau-t})\|^2
+\beta\suo\int_{-t}^0\suo\me^{\lam s
+2k\int^0_sz(\theta_r\w)\di r-2kz(\theta_s\w)}\|
\me^{kz(\theta_{s}\w)}v(s\!+\!\tau)\!-\!k^{-1}h(s\!+\!\tau)\|_q^q\di s\notag\\
\leqslant&\cR^2(\lam,\tau,\w),
\label{3.29z}\end{align}
where $\cR(\lam',\tau,\w)$ is defined for each $\lam'\in[\lam_0,\lam]$ as
\be\label{3.30z}
\cR^2(\lam',\tau,\w):=1+c\int_{-\8}^0\me^{\lam' s+2k\int_s^0z(\theta_r\w)\di r-2kz(\theta_s\w)}
\cQ_1(s+\tau)\di s<+\8,
\ee
and $\cQ_1(s)$ is given by \eqref{3.19z}.
Moreover, for each $(\tau,\w)\in\R\X\W$, let
\be\label{3.31z}K(\tau,\w):=\{u\in\ell^2:\|u\|\leqslant\me^{kz(\w)}
\cR(\lam,\tau,\w)+k^{-1}\|h(\tau)\|\}.\ee
Then
\be\label{3.32z}
\lim_{s\ra-\8}\me^{\lam_0 s+2k\int^{0}_{s}z(\theta_r\w)\di r
-2kz(\theta_s\w)}\|K(\tau+s,\theta_s\w)\|^2=0.
\ee
\el
\bo Replacing $\tau$, $\sig$ and $\w$ in \eqref{3.20z} by $\tau-t$, $\tau-t$
and $\theta_{-\tau}\w$, respectively, we obtain by variable substitution that,
\begin{align}&
\|v(s,\tau-t,\theta_{-\tau}\w,v_{\tau-t})\|^2\notag\\
&+\beta\int_{-t}^{s-\tau}\me^{\lam(\vsig+\tau-s)+2k\int_\vsig^{s-\tau}z(\theta_{r}\w)\di r
-2kz(\theta_{\vsig}\w)}\|\me^{kz(\theta_{\vsig}\w)}v(\vsig+\tau)-k^{-1}h(\vsig+\tau)\|_q^q
\di\vsig\notag\\
\leqslant&\me^{\lam(\tau-t-s)+2k\int_{-t}^{s-\tau}z(\theta_{r}\w)\di r}\|v_{\tau-t}\|^2\notag\\
&+c\int_{-t}^{s-\tau}\me^{\lam(\vsig+\tau-s)+2k\int_\vsig^{s-\tau}z(\theta_r\w)\di r
-2kz(\theta_\vsig\w)}\cQ_1(\vsig+\tau)\di\vsig.\label{3.33z}
\end{align}
When $s=\tau$ and $\lam'\in[\lam_0,\lam]$, it follows from the inequality \eqref{3.33z} that
\begin{align}&
\|v(\tau,\tau-t,\theta_{-\tau}\w,v_{\tau-t})\|^2\!+\!\beta\suo\int_{-t}^{0}\suo\me^{\lam s
+2k\int_s^0z(\theta_{r}\w)\di r
-2kz(\theta_{s}\w)}\|\me^{kz(\theta_{s}\w)}v(s+\tau)-k^{-1}h(s+\tau)\|_q^q
\di s\notag\\
\leqslant&\me^{-\lam' t+2k\int_{-t}^0z(\theta_{r}\w)\di r}\|v_{\tau-t}\|^2
+c\int_{-t}^0\me^{\lam' s+2k\int_s^0z(\theta_r\w)\di r-2kz(\theta_s\w)}
\cQ_1(s+\tau)\di s.\label{3.34z}
\end{align}
First we note by \eqref{2.6n} that there is $s_0:=s_0(\lam')<0$ such that when $s<s_0$,
\be\label{3.35z}
\min\left\{2kz(\theta_s\w),2k\int_{0}^sz(\theta_{r}\w)\di r\right\}>\frac{(\lam'-\lam_1)s}{2},\ee
and then it obviously holds that
\be\label{3.36z}\lam' s+2k\int_s^0z(\theta_r\w)\di r-2kz(\theta_s\w)
<\lam' s-(\lam'-\lam_1)s=\lam_1s.\ee
Thus by the selection $\me^{kz(\theta_{-t}\w)}v_{\tau-t}-k^{-1}h(\tau-t)\in D(\tau-t,\theta_{-t}\w)$
with $D\in\cD$, we see that when $-t<s_0$, by \eqref{3.27z} and \eqref{3.28z},
\begin{align}\label{3.37z}&\me^{-\lam' t+2k\int_{-t}^0z(\theta_{r}\w)\di r}\|v_{\tau-t}\|^2\notag\\
\leqslant&2\me^{-\lam' t+2k\int_{-t}^0z(\theta_{r}\w)\di r-2kz(\theta_{-t}\w)}
\(\|\me^{kz(\theta_{-t}\w)}v_{\tau-t}-k^{-1}h(\tau-t)\|^2+k^{-2}\|h(\tau-t)\|^2\)\notag\\
\leqslant&2\me^{-\lam_0t+2k\int_{-t}^0z(\theta_{r}\w)\di r
-2kz{(\theta_{-t}\w)}}\|D(\tau-t,\theta_{-t}\w)\|^2
+2k^{-2}\me^{-\lam_1t}\|h(\tau-t)\|^2\ra0,\end{align}
as $t\ra+\8$.
Consequently \eqref{3.29z} follows immediately and
the boundedness in \eqref{3.30z} follows from the hypotheses \textbf{(G2)}, \textbf{(H2)}
and \eqref{3.36z}.

Next we show \eqref{3.32z}.
By \eqref{3.35z} and \eqref{3.36z}, we know that when $s<s_0(\lam_0)$,
\be\label{3.38z}\lam_0 s+2k\int_s^0z(\theta_r\w)\di r-2kz(\theta_s\w)<\lam_1s
\hs\mb{and}\hs\lam_0 s+2k\int_s^0z(\theta_r\w)\di r<\lam_1s.\ee
Then by \eqref{3.27z} and \eqref{3.30z}, it yields that
\begin{align*}&\lim_{s\ra-\8}
\me^{\lam_0 s+2k\int^{0}_{s}z(\theta_r\w)\di r-2kz(\theta_s\w)}\|K(\tau+s,\theta_s\w)\|^2\\
\leqslant&2\lim_{s\ra-\8}\me^{\lam_0 s+2k\int^{0}_{s}z(\theta_r\w)\di r-2kz(\theta_s\w)}
\(\me^{2 kz(\theta_s\w)}\cR^2(\lam,\tau+s,\theta_s\w)+k^{-2}\|h(\tau+s)\|^2\)\\
\leqslant&c\lim_{s\ra-\8}\(\me^{\lam_1s}+\me^{\lam_0 s+2k\int^{0}_{s}z(\theta_r\w)\di r}
\int_{-\8}^0\me^{\lam\vsig+2k\int_\vsig^0z(\theta_{r+s}\w)\di r
-2kz(\theta_{\vsig+s}\w)}\cQ_1(\vsig+\tau+s)\di\vsig\)\\
&+2k^{-2}\lim_{s\ra-\8}\me^{\lam_1s}\|h(\tau+s)\|^2\\
=&c\lim_{s\ra-\8}\me^{\lam_0s+2k\int_s^0z(\theta_r\w)\di r}
\int_{-\8}^s\me^{\lam_0(\vsig-s)+2k\int_\vsig^sz(\theta_r\w)\di r
-2kz(\theta_\vsig\w)}\cQ_1(\vsig+\tau)\di\vsig\\
=&c\lim_{s\ra-\8}\int_{-\8}^s\me^{\lam_0\vsig+2k\int_\vsig^0z(\theta_r\w)\di r
-2kz(\theta_\vsig\w)}\cQ_1(\vsig+\tau)\di\vsig=0,
\end{align*}
by \eqref{3.38z}, \eqref{3.25z} and \eqref{3.26z} and proves \eqref{3.32z}.
The proof is finished.
\eo

Next we show the pullback $\cD$-asymptotic compactness of $\vp$.
Actually, in the case for lattice systems on $\ell^q$, $q\geqslant1$,
we only need to show that $\vp$ is \emph{pullback $\cD$-asymptotically null}
(see \cite{HSZ11,WJ24}),
i.e., for each $\ve>0$ and $D\in\cD$, there exist $T(\ve,\tau,\w,D)>0$ and $N(\ve,\tau,\w,D)\in\N$
such that for all $t\geqslant T(\ve,\tau,\w,D)$,
$$\sup_{u\in D(\tau,\w)}\sum_{|i|>N(\ve,\tau,\w,D)}|\vp_i(t,\tau-t,\theta_{-t}\w,u)|^q<\ve.$$
By the discussions in \cite{HSZ11},
one can infer that the pullback $\cD$-asymptotic nullity of $\vp$
ensures the pullback asymptotic compactness of $\vp$.

To prove the pullback $\cD$-asymptotic nullity of $\vp$ in $\ell^2$,
we need to introduce a smooth non-decreasing cut-off function $\rho:\R^+\ra[0,1]$ such that
$$\rho(s)=\left\{\ba{ll}0,&\mb{ for }s\in[0,1],\\1,&\mb{ for }s\in[2,+\8),
\ea\right.\hs\mb{and the derivative }\rho'(s)\leqslant c_0,$$
for all $s\in\R^+$ and some constant $c_0>0$.
For notational brevity, define for $N\in\N^+$,
\be\label{3.38y}\rho_i^N=\rho\(\frac{|i|}{N}\)\hs\mb{and}\hs\rho^N=\(\rho_i^N\)_{i\in\Z}.\ee
The following lemma guarantees the pullback $\cD$-asymptotic nullity of $\vp$ in $\ell^2$.

\bl\label{le3.5} Let the hypotheses \textbf{(F)}, \textbf{(G1)}, \textbf{(G2)},
\textbf{(H1)} and \textbf{(H2)} hold and
$\me^{kz(\theta_{-t}\w)}v_{\tau-t}-k^{-1}h(\tau-t)\in D(\tau-t,\theta_{-t}\w)$ with $D\in\cD$.
Then for every $\ve>0$ and $\bP$-a.s. $\w\in\W$, there exist $T_1=T_1(\ve,\tau,\w,D)>0$
and $I=I(\ve,\tau,\w,D)\in\N$ such that when $t>T_1$, it holds that
\be\sum_{|i|>I}|
v_i(\tau,\tau-t,\theta_{-\tau}\w,v_{\tau-t})|^2<\ve.
\label{3.39z}\ee
\el
\bo Consider the solution $v(t,\tau,\w,v_\tau)$ of \eqref{3.8z}.
Let $N$ be a positive integer with $N\geqslant3$ and $w(t):=\rho^N\otimes v(t)$.
We only discuss $\|w(\tau,\tau-t,\theta_{-\tau}\w,v_{\tau-t})\|^2$ on account of the evident fact
$$\|w(\tau,\tau-t,\theta_{-\tau}\w,v_{\tau-t})\|^2\geqslant
\sum_{|i|>2N}|v_i(\tau,\tau-t,\theta_{-\tau}\w,v_{\tau-t})|^2.$$
In the following, we give some necessary estimates for $(F(t,\w,v),w)$ in advance.

First by Cauchy inequality and \eqref{2.3n}, we have
\begin{align}(A^2v,w)=&(Av,Aw)=\sum_{i\in\Z}(2v_i-v_{i+1}-v_{i-1})
(2\rho_i^Nv_i-\rho_{i+1}^Nv_{i+1}-\rho_{i-1}^Nv_{i-1})\notag\\
=&\sum_{i\in\Z}\rho_i^N(Av)_i^2+\sum_{i\in\Z}(\rho_i^N-\rho_{i+1}^N)(Av)_iv_{i+1}
+\sum_{i\in\Z}(\rho_i^N-\rho_{i-1}^N)(Av)_iv_{i-1}\notag\\
\geqslant&-\frac{c_0}{N}\(\sum_{i\in\Z}|(Av)_i|(|v_{i+1}|+|v_{i-1}|)\)\notag\\
\geqslant&-\frac{c_0}{N}\(\frac14\|Av\|^2+4\|v\|^2\)\geqslant-\frac{8c_0}{N}\|v\|^2,
\label{3.39y}\end{align}
and hence
$$
(F_1(t,\w,v),w)\leqslant\frac{8c_0}{N}\|v\|^2+(kz(\theta_t\w)-\lam)\|w\|^2.
$$
For $F_2(t,\w,v)$, it is obvious that
\begin{align*}\(F_2(t,\w,v),w\)=&\(F_2(t,\w,v),\rho^N\otimes(v-k^{-1}\me^{-kz(\theta_t\w)}h(t))
+k^{-1}\me^{-kz(\theta_t\w)}\rho^N\otimes h(t)\)\\
=&\me^{-2kz(\theta_t\w)}\(-\rho^N\otimes Af(\me^{kz(\theta_t\w)}v-k^{-1}h(t)),
\rho^N\otimes(\me^{kz(\theta_t\w)}v-k^{-1}h(t))\)\\
&+\me^{-2kz(\theta_t\w)}\(g(t,\me^{kz(\theta_t\w)}v-k^{-1}h(t)),
\rho^N\otimes(\me^{kz(\theta_t\w)}v-k^{-1}h(t))\)\\
&+k^{-1}\me^{-2kz(\theta_t\w)}\(-\rho^N\otimes Af(\me^{kz(\theta_t\w)}v-k^{-1}h(t)),
\rho^N\otimes h(t))\)\\
&+k^{-1}\me^{-2kz(\theta_t\w)}\(g(t,\me^{kz(\theta_t\w)}v-k^{-1}h(t)),
\rho^N\otimes h(t)\)\notag\\
:=&J_1+J_2+J_3+J_4.\end{align*}
Below we note that when $|i|<N$, $\rho_{i\pm1}^N=\rho_i^{N-1}=0$; when $|i|\geqslant N$,
$$\frac{|i|}{N-1}-\frac{|i\pm 1|}{N}=\left\{\ba{ll}\disp\frac{i\pm(N-1)}{N(N-1)},&i>0,\\[2ex]
\disp\frac{-i\mp(N-1)}{N(N-1)},&i<0\ea\right.>0.$$
Recalling the non-decreasing of $\rho$, we deduce that
$$\max\{\rho_i^N,\rho_{i-1}^N,\rho^{N}_{i+1}\}\leqslant\rho^{N-1}_i.$$
Hence for each $x\in\ell^q$ with each $q\geqslant1$,
\begin{align}\|\rho^N\otimes Ax\|_q^q=&\sum_{i\in\Z}\rho^N_i|2x_i-x_{i+1}-x_{i-1}|^q\notag\\
\leqslant&2^{2q-2}\sum_{i\in\Z}\(2\rho^N_i+\rho^N_{i-1}+\rho^N_{i+1}\)|x_i|^q\notag\\
\leqslant&2^{2q}\sum_{|i|\geqslant N}\rho^{N-1}_i|x_i|^q
=4^{q}\|\rho^{N-1}\otimes x\|_q^q.\label{3.40z}
\end{align}
Then by the hypotheses \textbf{(F)}, \textbf{(G1)}, \textbf{(H1)},
\eqref{2.2n} and Young's inequality,
we know that
\begin{align*}\me^{2kz(\theta_t\w)}|J_1|\leqslant&
\|\rho^N\otimes Af(\me^{kz(\theta_t\w)}v-k^{-1}h(t))\|_1
\|\me^{kz(\theta_t\w)}v-k^{-1}h(t)\|_q\\
\leqslant&4\|\rho^{N-1}\otimes f(\me^{kz(\theta_t\w)}v-k^{-1}h(t))\|_1
\|\me^{kz(\theta_t\w)}v-k^{-1}h(t)\|_q\\
\leqslant&4\|\me^{kz(\theta_t\w)}v-k^{-1}h(t)\|_q
\sum_{i\in\Z}\rho^{N-1}_i
\(\alpha_{1i}|\me^{kz(\theta_t\w)}v_i-k^{-1}h_i(t)|^{q-1-\epsilon}+\alpha_{2i}\)\\
\leqslant&4\|\rho^{N-1}\!\!\otimes\alpha_{1}\|_1\|\me^{kz(\theta_t\w)}v-k^{-1}h(t)\|_q^{q-\epsilon}
\!\!+\!4\|\rho^{N-1}\!\!\otimes\alpha_{2}\|_1\|\me^{kz(\theta_t\w)}v-k^{-1}h(t)\|_q\\
\leqslant&c\(1+\|\me^{kz(\theta_t\w)}v-k^{-1}h(t)\|_q^q\)
\(\|\rho^{N-1}\otimes\alpha_{1}\|_1+\|\rho^{N-1}\otimes\alpha_{2}\|_1\),\\
\me^{2kz(\theta_t\w)}J_2\leqslant&
\sum_{i\in\Z}\rho_i^N(-\beta|\me^{kz(\theta_t\w)}v_i-k^{-1}h_i(t)|^q+\psi_{1i}(t))\\
\leqslant&-\beta\|\rho^N\otimes(\me^{kz(\theta_t\w)}v-k^{-1}h(t))\|_q^q
+\|\rho^N\otimes\psi_1(t)\|_1,\\
\me^{2kz(\theta_t\w)}|J_3|\leqslant&
c\|\rho^{N-1}\otimes f(\me^{kz(\theta_t\w)}v-k^{-1}h(t))\|_1\|\rho^N\otimes h(t)\|_q\\
\leqslant&c\|h(t)\|_q\(\|\rho^{N-1}\otimes\alpha_{1}\|_1\|
\me^{kz(\theta_t\w)}v-k^{-1}h(t)\|_q^{q-1-\epsilon}
+\|\rho^{N-1}\otimes\alpha_2\|_1\)\notag\\
\leqslant&c\|\rho^{N-1}\otimes\alpha_{1}\|_1\(1+\|h(t)\|_q^q
+\|\me^{kz(\theta_t\w)}v-k^{-1}h(t)\|_q^q\)\\
&+c\|\rho^{N-1}\otimes\alpha_2\|_1\(1+\|h(t)\|_q^q\)\hs\mb{and}\\
\me^{2kz(\theta_t\w)}|J_4|\leqslant&c\sum_{i\in\Z}
\(\psi_{3i}(t)|\me^{kz(\theta_t\w)}v_i-k^{-1}h_i(t)|^{q-1}+\psi_{4i}(t)\)\rho^N_i|h_i(t)|\\
\leqslant&\beta\|\rho^N\otimes(\me^{kz(\theta_t\w)}v-k^{-1}h(t))\|_q^q
+c\|\rho^N\otimes\psi_4(t)\|^2\\
&+c\|\psi_3(t)\|^q_{\8}\|\rho^N\otimes h(t)\|_q^q+c\|\rho^N\otimes h(t)\|^2.
\end{align*}
As to $F_3(t,\w)$, we can easily see by \eqref{3.40z} that
\begin{align*}|(F_3(t,\w),w)|\leqslant&c\me^{-2kz(\theta_t\w)}\(\|\rho^N\otimes A^2h(t)\|^2
+\|\rho^N\otimes h(t)\|^2+\|\rho^N\otimes h'(t)\|^2\)+\frac{\lam}{2}\|w\|^2\\
\leqslant&c\me^{-2kz(\theta_t\w)}\(\|\rho^{N-2}\otimes h(t)\|^2+\|\rho^N\otimes h(t)\|^2
+\|\rho^N\otimes h'(t)\|^2\)+\frac{\lam}{2}\|w\|^2\\
\leqslant&c\me^{-2kz(\theta_t\w)}\(\|\rho^{N-2}\otimes h(t)\|^2
+\|\rho^N\otimes h'(t)\|^2\)+\frac{\lam}{2}\|w\|^2.
\end{align*}

Similar to the setting \eqref{3.19z}, we let
\begin{align*}\cQ^N(t,\w,v):=&\(\|\rho^{N-1}\otimes\al_1\|_1+\|\rho^{N-1}\otimes\al_2\|_1\)
\(1+\|h(t)\|_q^q+\|\me^{kz(\theta_t\w)}v-k^{-1}h(t)\|_q^q\)\notag\\
&+\|\rho^N\otimes\psi_1(t)\|_1+\|\rho^N\otimes\psi_4(t)\|^2+\|\psi_3(t)\|^q_{\8}
\|\rho^N\otimes h(t)\|_q^q\\
&+\|\rho^{N-2}\otimes h(t)\|^2+\|\rho^N\otimes h'(t)\|^2.
\end{align*}
Then by taking the inner product of \eqref{3.8z} and $w$ and using the estimates above,
we obtain
\be\label{3.41z}\frac{\di\|w(t)\|^2}{\di t}+(\lam-2kz(\theta_t\w))\|w(t)\|^2\leqslant
\frac{16c_0}{N}\|v(t)\|^2+c\me^{-2kz(\theta_t\w)}\cQ^N(t,\w,v(t)).\ee

Applying Gronwall's lemma to \eqref{3.41z} over $[\sig,s]$ and
substituting $\tau$, $\sig$, $s$ and $\w$ by $\tau-t$, $\tau-t$, $\tau$
and $\theta_{-\tau}\w$, respectively,
we can get
\begin{align*}
&\|w(\tau,\tau-t,\theta_{-\tau}\w,v_{\tau-t})\|^2\\
\leqslant&\me^{-\lam t+2k\int_{-t}^0 z(\theta_{r}\w)\di r}\|v_{\tau-t}\|^2
+\frac{16c_0}{N}\int_{-t}^0\me^{\lam\vsig+2k\int_{\vsig}^0z(\theta_{r}\w)\di r}
\|v(\vsig+\tau)\|^2\di\vsig\\
&+c\int_{-t}^0\me^{\lam\vsig+2k\int_{\vsig}^0z(\theta_{r}\w)\di r-2kz(\theta_{\vsig}\w)}
\cQ^N(\vsig+\tau,\theta_{-\tau}\w,v(\vsig+\tau))\di \vsig\\
:=&L_1+L_2+L_3.\end{align*}
The conclusion \eqref{3.37z} has assured that $L_1$ vanishes as $t\ra+\8$.
For $L_2$, we observe by the discussion from \eqref{3.34z} to \eqref{3.37z} that
there is $T_{\lam_0}$ large enough such that when $t>T_{\lam_0}$,
\begin{align*}&\|v(\vsig+\tau,\tau-t,\theta_{-\tau}\w,v_{\tau-t})\|^2\\
\leqslant&\me^{-\lam(t+\vsig)+2k\int_{-t}^{\vsig}z(\theta_{r}\w)\di r}\|v_{\tau-t}\|^2
+c\int_{-t}^{\vsig}\me^{\lam(\sig-\vsig)
+2k\int_\sig^{\vsig}z(\theta_r\w)\di r-2kz(\theta_\sig\w)}\cQ_1(\sig+\tau)\di\sig\\
\leqslant&\me^{-\lam_0\vsig}\(\me^{-\lam_0t
+2k\int_{-t}^{\vsig}z(\theta_{r}\w)\di r}\|v_{\tau-t}\|^2
+c\int_{-\8}^0\me^{\lam_0\sig+2k\int_\sig^{\vsig}z(\theta_r\w)\di r-2kz(\theta_\sig\w)}
\cQ_1(\sig+\tau)\di\sig\)\\
\leqslant&\me^{-\lam_0\vsig-2k\int_\vsig^0z(\theta_r\w)\di r}\cR^2(\lam_0,\tau,\w)
\end{align*}
and hence
$$L_2\leqslant\frac{16c_0}{N}\cR^2(\lam_0,\tau,\w)
\int_{-\8}^0\me^{(\lam-\lam_0)\vsig}\di\vsig\ra0
$$
as $N\ra+\8$.
At last, we consider $L_3$.
Indeed, based on the hypotheses \textbf{(F)}, \textbf{(G1)}, \textbf{(G2)},
\textbf{(H1)} and \textbf{(H2)},
what is essentially worth attention in $L_3$ lies in the term
$$\(\|\rho^{N-1}\otimes\al_1\|_1+\|\rho^{N-1}\otimes\al_2\|_1\)\suo
\int_{-t}^0\me^{\lam\vsig+2k\int_\vsig^0z(\theta_r\w)\di r-2kz(\theta_\vsig\w)}
\|\me^{kz(\theta_t\w)}v(\vsig+\tau)-k^{-1}h(t)\|_q^q\di\vsig,$$
in which, however, the integral is uniformly bounded by $\cR^2(\lam,\tau,\w)/\beta$
when $t>T_0$ (given in Lemma \ref{le3.4}).
This is sufficient to guarantee the term mentioned above to vanish when $N\ra+\8$.
Because of the boundedness given in \eqref{3.30z},
it is easy to see that each of other terms of $L_3$ tends to zero as $N\ra+\8$;
for example, as to the term including $\|\rho^N\otimes\psi_1(\vsig+\tau)\|_1$,
we know by \eqref{3.25z} and \eqref{3.36z} that
\be\label{3.42z}\int_{-t}^0\me^{\lam\vsig+2k\int_\vsig^0z(\theta_r\w)\di r-2kz(\theta_\vsig\w)}
\|\psi_1(\vsig+\tau)\|_1\di\vsig<+\8,\ee
and when $\psi_1$ is replaced by $\rho^N\otimes\psi_1$ in \eqref{3.42z},
the left side of \eqref{3.42z} becomes less than its remainder term of
the corresponding series (using the norm of $\psi_1(\vsig+\tau)$),
which surely tends to $0$ as $N\ra+\8$.

All in all, we have effectively proved that for every $\ve>0$ and $\bP$-a.s. $\w\in\W$,
there exist $T_1=T_1(\ve,\tau,\w,D)>0$ and $N=N(\ve,\tau,\w,D)\in\N$ such that when $t>T_1$,
it holds that $\|\tau,\tau-t,\theta_{-\tau}\w,v_{\tau-t}\|^2<\ve$,
which is exactly the conclusion we desire.
The proof is finished.
\eo

Now we are well prepared to investigate the existence of pullback random
$\cD$-attractors in $\ell^2$.
\bt\label{th3.6} Let the hypotheses \textbf{(F)}, \textbf{(G1)}, \textbf{(G2)},
\textbf{(H1)} and \textbf{(H2)} hold.
Then $\vp$ has a unique pullback random $\cD$-attractor $\cA$ and $\cA\in\cD$.
\et
\bo We only need to check the two conditions presented in Theorem \ref{th2.4n}.

Let $K:=\{K(\tau,\w)\}_{(\tau,\w)\in\R\X\W}$ be given in \eqref{3.31z}.
The randomness of $K$ follows from Proposition \ref{pro2.2n}.
By \eqref{3.29z} and \eqref{3.23z},
it can be inferred that for each $u_{\tau-t}\in D(\tau-t,\theta_{-t}\w)$ with $D\in\cD$,
$\tau\in\R$ and $\bP$-a.s. $\w\in\W$,
\begin{align*}\|u(\tau,\tau-t,\theta_{-\tau}\w,u_{\tau-t})\|=&
\|\me^{kz(\w)}v(\tau,\tau-t,\theta_{-\tau}\w,v_{\tau-t})-k^{-1}h(\tau)\|\\
\leqslant&\me^{kz(\w)}\cR(\lam,\tau,\w)+k^{-1}\|h(\tau)\|,\end{align*}
which indicates that $K$ is a pullback $\cD$-absorbing set for $\vp$ in $\ell^2$.
By \eqref{3.28z} and \eqref{3.32z}, we know $K\in\cD$.

To show the pullback $\cD$-asymptotic compactness of $\vp$, we adopt Lemma \ref{le3.5}.
Consider the initial datum $u_{\tau-t}\in D(\tau-t,\theta_{-t}\w)$ for each $D\in\cD$.
Then by \eqref{3.24z} and setting $v_{\tau-t}=\me^{-kz(\theta_{-t}\w)}(u_{\tau-t}+k^{-1}h(\tau-t))$,
we have
\begin{align*}&\sum_{|i|>N}\left|\(\vp(t,\tau-t,\theta_{-t}\w,u_{\tau-t})\)_i\right|^2\\
=&\sum_{|i|>N}\left|\me^{kz(\w)}v_i(\tau,\tau-t,
\theta_{-\tau}\w,v_{\tau-t})-k^{-1}h_i(\tau)\right|^2\\
\leqslant&2\(\me^{2kz(\w)}\sum_{|i|>N}\left|v_i(\tau,\tau-t,\theta_{-\tau}\w,v_{\tau-t})\right|^2
+k^{-2}\sum_{|i|>N}|h_i(\tau)|^2\).
\end{align*}
Recall \eqref{3.39z} and $h(\tau)\in\ell^2$ and we can directly obtain
the pullback $\cD$-asymptotic nullity
of $\vp$ and hence the pullback $\cD$-asymptotic compactness.

Applying Theorem \ref{th2.4n}, we then obtain a unique pullback random
$\cD$-attractor $\cA\in\cD$ for $\vp$ in $\ell^2$.
The proof is thus finished.
\eo

\subsection{Existence of invariant sample measures}\label{ss3.4}

In this subsection, we investigate the existence of invariant sample measures
for $\vp$ in $\ell^2$.
According to Proposition \ref{pro2.6n} and the discussion in Section 3,
we are only necessary to study the boundedness of the mapping
$\tau\mapsto\vp(t-\tau,\tau,\theta_{\tau}\w,\phi)$ over $(-\8,t]$
and the joint continuity of the mapping $(\tau,\phi)\mapsto\vp(t-\tau,\tau,\theta_\tau\w,\phi)$
in $\ell^2$ for fixed $t\in\R$ and $\bP$-a.s. $\w\in\W$.

As before, we start our discussion still from $v(t,\tau,\w,\phi)$.
\bl\label{le3.7} Let the hypotheses \textbf{(F)}, \textbf{(G1)}, \textbf{(G2)}, \textbf{(H1)}
and \textbf{(H2)} hold.
Then for each fixed $(t,\phi)\in\R\X\ell^2$, $\bP$-a.s., $\w\in\W$ and all $\tau\leqslant t$,
\be\label{3.43z}\|v(t,\tau,\w,\phi)\|\leqslant\cM(t,\w,\phi)<+\8,\ee
where
\begin{align}\cM^2(t,\w,\phi):=&\|\phi\|^2\max_{\tau\leqslant t}
\me^{\lam(\tau-t)+2k\int_\tau^tz(\theta_r\w)\di r}\notag\\
&+c\int_{-\8}^t\me^{\lam(\vsig-t)+2k\int_{\vsig}^tz(\theta_r\w)\di r
-2kz(\theta_\vsig\w)}\cQ_1(\vsig)\di\vsig.
\label{3.44z}\end{align}
\el
\bo Choosing $\sig=\tau$, $s=t$ and $v_\tau=\phi$ in \eqref{3.20z},
we have
\be\|v(t)\|^2\leqslant\me^{\lam(\tau-t)+2k\int_\tau^tz(\theta_r\w)\di r}\|\phi\|^2
+c\int_\tau^t\me^{\lam(\vsig-t)+2k\int_{\vsig}^tz(\theta_r\w)\di r
-2kz(\theta_\vsig\w)}\cQ_1(\vsig)\di\vsig,
\label{3.45z}\ee
which is uniformly dominated by $\cM^2(t,\w,\phi)$ for each $\tau\leqslant t$.
Thanks to \eqref{3.35z}, \eqref{3.36z} and all the hypotheses given above,
we see that
$$\lim_{\tau\ra-\8}\me^{\lam(\tau-t)+2k\int_\tau^tz(\theta_r\w)\di r}=0$$
and the integral term in \eqref{3.44z} is bounded.
Consequently, $\cM^2(t,\w,\phi)<+\8$ and the proof is accomplished.
\eo
\bl\label{le3.8} Let the hypotheses \textbf{(F)}, \textbf{(G1)}, \textbf{(G2)},
\textbf{(H1)} and \textbf{(H2)} hold.
Fix $(t,\phi)\in\R\X\ell^2$ and $\bP$-a.s., $\w\in\W$.
Then for every $\ve>0$, there exists $\de_1=\de_1(\ve,t,\w,\phi)>0$ such that when
$\tau\in(t-\de_1,t)$,
\be\label{3.46z}\|v(t,\tau,\w,\phi)-\phi\|<\ve.\ee
\el
\bo First, we notice from \eqref{3.8z} and by H\"older inequality that
\be\label{3.47z}\|v(t,\tau,\w,\phi)-\phi\|^2=\left\|\int_{\tau}^{t}F(s,\w,v(s))\di s\right\|^2
\leqslant(t-\tau)\int_{\tau}^t\|F(s,\w,v(s))\|^2\di s.\ee
Therefore, we need to estimate $\|F(s,\w,v(s))\|$.
Note by Young's inequality and embeddings \eqref{2.1} that
\be\|F_1(s,\w,v(s))\|^2\leqslant c(1+|z(\theta_s\w)|^2)\|v(s)\|^2,\ee
\begin{align}\me^{2kz(\theta_s\w)}\|F_2(s,\w,v(s))\|^2
\leqslant&
c\|f(\me^{kz(\theta_s\w)}v(s)-k^{-1}h(s))\|^2+\|g(t,\me^{z(\theta_s\w)}v(s)-k^{-1}h(s))\|^2\notag\\
\leqslant&c\sum_{i\in\Z}\(\alpha_{1i}^2|\me^{kz(\theta_s\w)}v_i-k^{-1}h_i(s)|^{2q-2-2\epsilon}
+\alpha_{2i}^2\)\notag\\
&+c\sum_{i\in\Z}\(|\psi_{3i}(s)|^2|\me^{kz(\theta_s\w)}v_i-k^{-1}h_i(s)|^{2q-2}
+|\psi_{4i}(s)|^2\)\notag\\
\leqslant&c\(\|\al_1\|_1^{\frac{2q}{1+\epsilon}}+\|\me^{kz(\theta_s\w)}v-k^{-1}h(s)\|^{2q}
+\|\al_2\|_1^2\)\notag\\
&+c\(\|\psi_3(s)\|_\8^{2q}+\|\me^{kz(\theta_s\w)}v-k^{-1}h(s)\|^{2q}+\|\psi_4(s)\|^2\)\notag\\
\leqslant&c\(\|\al_1\|_1^{\frac{2q}{1+\epsilon}}+\|\al_2\|_1^2+\|\psi_3(s)\|_\8^{2q}
+\|\psi_4(s)\|^2+\|h(s)\|^{2q}\)\notag\\
&+c\me^{2qkz(\theta_s\w)}\|v(s)\|^{2q}\hs\mb{and}
\end{align}
\be\me^{2kz(\theta_s\w)}\|F_3(s,\w)\|^2\leqslant c\(\|h(s)\|^2+\|h'(s)\|^2\).\ee

By \eqref{3.43z}, we have, for $s\in[\tau,t]$,
\begin{align}&\|v(s)\|^2\leqslant\cM^2(s,\w,\phi)\notag\\
=&\|\phi\|^2\max_{\tau\leqslant s}\me^{\lam(\tau-s)+2k\int_\tau^sz(\theta_r\w)\di r}
+c\int_{-\8}^s\me^{\lam(\vsig-s)+2k\int_{\vsig}^sz(\theta_r\w)\di r-2kz(\theta_\vsig\w)}
\cQ_1(\vsig)\di\vsig\notag\\
\leqslant&\me^{\lam(t-s)-2k\int_s^tz(\theta_r\w)\di r}\cM^2(t,\w,\phi),\label{3.51z}
\end{align}
Setting
$$\cQ_2(s):=\|\al_1\|_1^{\frac{2q}{1+\epsilon}}+\|\al_2\|_1^2+\|\psi_3(s)\|_\8^{2q}
+\|\psi_4(s)\|^2+\|h(s)\|^{2q}+\|h(s)\|^2+\|h'(s)\|^2,$$
we deduce by the inequalities from \eqref{3.47z}-\eqref{3.51z} that when $t>\tau>t-1$,
\begin{align}&\frac{1}{t-\tau}\|v(t,\tau,\w,\phi)-\phi\|^2\notag\\
\leqslant&c\int_{\tau}^t\(\(1+|z(\theta_s\w)|^2\)\|v(s)\|^2
+\me^{2(q-1)kz(\theta_s\w)}\|v(s)\|^{2q}
+\me^{-2kz(\theta_s\w)}\cQ_2(s)\)\di s\notag\\
\leqslant&c\cM^2(t,\w,\phi)\int_{t-1}^t\(1+|z(\theta_s\w)|^2\)
\me^{\lam(t-s)-2k\int_s^tz(\theta_r\w)\di r}\di s
+c\int_{t-1}^t\me^{-2kz(\theta_s\w)}\cQ_2(s)\di s\notag\\
&+c\cM^{2q}(t,\w,\phi)\int_{t-1}^t\me^{\lam q(t-s)
-2qk\int_s^tz(\theta_r\w)\di r+2(q-1)kz(\theta_s\w)}\di s.
\label{3.52z}\end{align}
The boundedness of the right side of \eqref{3.52z} is obvious except the integral below included
in the second term
$$\int_\tau^t\me^{-2kz(\theta_s\w)}\|h(s)\|^{2q}\di s,$$
which is also bounded, since $h:\R\ra\ell^2$ is derivable and whereby continuous,
and then $\|h(s)\|^{2q}$ is continuous in $s$.

Consequently, for all $\ve>0$,
the boundedness of \eqref{3.52z} enables us to pick a sufficiently small
$\de_1:=\de_1(\ve,t,\w,\phi)\in(0,1)$ such that when $\tau\in(t-\de_1,t)$,
$$\|v(t,\tau,\w,\phi)-\phi\|^2<\ve^2.$$
This is exactly \eqref{3.46z}. The proof is hence finished.
\eo
\bl\label{le3.9} Let the hypotheses \textbf{(F)}, \textbf{(G1)}, \textbf{(G2)},
\textbf{(H1)} and \textbf{(H2)} hold.
Then given each fixed $\tau,\,t\in\R$ with $t>\tau$, $\phi\in\ell^2$
and $\bP$-a.s. $\w\in\W$, for every $\ve>0$, there exists $\de_2=\de_2(\ve,t,\tau,\w,\phi)>0$
such that whenever $|\tau'-\tau|+\|\phi'-\phi\|<\de_2$,
\be\label{3.53z}\|v(t,\tau',\w,\phi')-v(t,\tau,\w,\phi)\|<\ve.\ee
As a result, the mapping $(\tau,\phi)\mapsto v(t,\tau,\w,\phi)$ is continuous
for each $t\geqslant\tau$ and $\bP$-a.s. $\w\in\W$.
\el

\bo Given $\tau,\,t\in\R$, $\phi\in\ell^2$ and $\w\in\W$ as stated,
we can use \eqref{3.21z} and conclude that
for all $\tau'\in(\tau_0,t)$ with $\tau_0=\tau-1$ and $\phi'$
with $\|\phi'\|\leqslant R_\phi:=\|\phi\|+1$,
$$\|v(t,\tau',\w,\phi')\|^2\leqslant\cP_2(\tau',t,\w,\|\phi'\|)
\leqslant\cP_2(\tau_0,t,\w,R_\phi).$$
Next we first consider the right continuity for \eqref{3.53z}, saying, the case when $\tau'>\tau$.
By \eqref{3.12z}, we know
\begin{align}\label{3.54z}&\frac{\di}{\di t}\|v(t,\tau',\w,\phi')-v(t,\tau,\w,\phi)\|^2\notag\\
\leqslant&2\cP_1\(\tau_0,t,\w,\sqrt{\cP_2(\tau_0,t,\w,R_\phi)}\)
\|v(t,\tau',\w,\phi')-v(t,\tau,\w,\phi)\|^2\notag\\
:=&2\~\cP_1(\tau_0,t,\w,R_\phi)\|v(t,\tau',\w,\phi')-v(t,\tau,\w,\phi)\|^2.\end{align}
Then it follows from \eqref{3.54z} that
\begin{align*}&\|v(t,\tau',\w,\phi')-v(t,\tau,\w,\phi)\|\leqslant
\|\phi'-v(\tau',\tau,\w,\phi)\|\me^{(t-\tau_0)\~\cP_1(\tau_0,t,\w,R_\phi)}\\
\leqslant&\(\|\phi'-\phi\|+\|\phi-v(\tau',\tau,\w,\phi)\|\)
\me^{(t-\tau_0)\~\cP_1(\tau_0,t,\w,R_\phi)}.
\end{align*}

Then the continuity of $v(t,\tau,\w,\phi)$ in $t$ given in Theorem \ref{th3.2} implies that
for each $\ve>0$, there is $\de'_2=\de'_2(\ve,t,\tau,\w,\phi)\in(0,1)$ such that
\eqref{3.53z} holds for all $(\tau',\phi')\in(\tau,+\8)\X\ell^2$
with $\tau'-\tau+\|\phi'-\phi\|<\de'_2$.
The left continuity is a similar conclusion by using Lemma \ref{le3.8}.
The proof is complete.
\eo

\bl\label{le3.10} Let the hypotheses \textbf{(F)}, \textbf{(G1)}, \textbf{(G2)},
\textbf{(H1)} and \textbf{(H2)} hold.
Then the mapping $\tau\mapsto u(t,\tau,\w,\phi)$ is bounded for $\tau\in(-\8,t]$,
and the mapping $(\tau,\phi)\ra u(t,\tau,\w,\phi)$ is continuous.\el
\bo Recalling \eqref{3.23z} and following Lemma \ref{le3.7},
we know that for all $\tau\in(-\8,t]$,
\begin{align*}&\|u(t,\tau,\w,\phi)\|=\|\me^{kz(\theta_t\w)}
v(t,\tau,\w,\me^{-kz(\theta_\tau\w)}\phi+k^{-1}h(\tau))-k^{-1}h(t)\|\\
\leqslant&\me^{kz(\theta_t\w)}\|v(t,\tau,\w,\me^{-kz(\theta_\tau\w)}\phi+k^{-1}h(\tau))\|
+k^{-1}\|h(t)\|.\end{align*}
By \eqref{3.45z}, we deduce that
\begin{align*}&\|v(t,\tau,\w,\me^{-kz(\theta_\tau\w)}\phi+k^{-1}h(\tau))\|^2\\
\leqslant&\me^{\lam(\tau-t)+2k\int_\tau^tz(\theta_r\w)\di r}
\|\me^{-kz(\theta_\tau\w)}\phi+k^{-1}h(\tau)\|^2
+c\int_\tau^t\me^{\lam(\vsig-t)+2k\int_{\vsig}^tz(\theta_r\w)\di r
-2kz(\theta_\vsig\w)}\cQ_1(\vsig)\di\vsig\\
\leqslant&c\me^{\lam(\tau-t)+2k\int_\tau^tz(\theta_r\w)\di r
-2kz(\theta_\tau\w)}\|\phi\|^2+c\me^{\lam(\tau-t)
+2k\int_\tau^tz(\theta_r\w)\di r}\|h(\tau)\|^2\\
&+c\int_\tau^t\me^{\lam(\vsig-t)+2k\int_{\vsig}^tz(\theta_r\w)\di r
-2kz(\theta_\vsig\w)}\cQ_1(\vsig)\di\vsig,
\end{align*}
which is uniformly bounded for all $\tau\in(-\8,t]$ by \eqref{3.27z} and \eqref{3.35z}.

The continuity of the mapping $(\tau,\phi)\ra u(t,\tau,\w,\phi)$ is
a direct result of the continuity of
$z(\theta_\tau\w)$ and $h(\tau)$ in $\tau$ and Lemma \ref{le3.9}.
The proof is complete.\eo

Based on Proposition \ref{pro2.6n}, Theorem \ref{th3.6} and Lemma \ref{le3.10}
and using the definition \eqref{3.24z},
we can conclude the existence of the invariant sample measures of $\vp$ on $\ell^2$ as follows.

\bt\label{th3.11}
Suppose that the hypotheses \textbf{(F)}, \textbf{(G1)}, \textbf{(G2)}, \textbf{(H1)}
and \textbf{(H2)} hold.
Let $\cA$ be the unique pullback random $\cD$-attractor of $\vp$
given in Theorem \ref{th3.6}.
Then for a given generalized Banach limit $\disp\LIM_{t\ra-\8}$ and a mapping $\xi:\R\X\W\ra\ell^2$
with $\{\xi(\tau,\w)\}_{(\tau,\w)\in\R\X\W}\in\cD$ such that
$\xi(\tau,\theta_\tau\w)$ is continuous in $\tau$,
there exists a family of Borel probability measures $\{\mu_{\tau,\w}\}_{(\tau,\w)\in\R\X\W}$
on $\ell^2$ such that the support of the measure $\mu_{\tau,\w}$ is contained in $\cA(\tau,\w)$ and
for all $\Upsilon\in\cC(\ell^2)$,
\begin{align}&
 \LIM_{\tau\ra-\8}\frac{1}{t-\tau}\int_{\tau}^{t}
 \Upsilon(\vp(t-s,s,\theta_s\w,\xi(s,\theta_s\w)))\di s
 =\int_{\cA(t,\theta_t\w)}\Upsilon(u)\mu_{t,\theta_t\w}(\di u)\notag\\
=&\int_{\ell^2}\Upsilon(u)\mu_{t,\theta_t\w}(\di u)
 =\LIM_{\tau\ra-\8}\frac{1}{t-\tau}\int_{\tau}^{t}\int_{\ell^2}\Upsilon(\vp(t-s,s,\theta_s\w,u))
 \mu_{s,\theta_s\w}(\di u)\di s.
\label{3.55z}\end{align}
Additionally, $\mu_{\tau,\w}$ is invariant in the sense that for all $\tau\in\R$ and
$\bP$-a.s. $\w\in\W$,
\be\label{3.56z}
 \int_{\cA(\tau+t,\theta_t\w)}\Upsilon(u)\mu_{\tau+t,\theta_t\w}(\di u)=
 \int_{\cA(\tau,\w)}\Upsilon(\vp(t,\tau,\w,u))\mu_{\tau,\w}(\di u),
 \hs\mb{for all }t\geqslant0.
\ee\et

The definition \eqref{3.24z} for $\vp$ also allows us to write \eqref{3.56z} with the solution $u$
in the following form,
\be\label{3.57z}
 \int_{\cA(\tau+t,\theta_t\w)}\Upsilon(u)\mu_{\tau+t,\theta_t\w}(\di u)=
 \int_{\cA(\tau,\w)}\Upsilon(u(t+\tau,\tau,\theta_{-\tau}\w,\phi))\mu_{\tau,\w}(\di\phi),
 \hs\mb{for all }t\geqslant0.
\ee

\subsection{Stochastic Liouville type theorem and sample statistical solutions}

In this subsection we prove the invariant sample measures given by Theorem \ref{th3.11}
satisfy the stochastic Liouville type theorem.

For writing simplification, we follow the notations in Subsection \ref{ss2.5} and let
\be\label{3.59z1}V:=\ell^2\cap\ell^q,\hs H=\ell^2\hs\mb{and}\hs V^*=\ell^2+\ell^p,\ee
where $q>2$ is given in the hypothesis \textbf{(F)} and $p\in[1,2)$ is the dual index of $q$,
i.e., $1/p+1/q=1$.
The norm of $u\in V$ is endowed as $\|u\|_V:=\|u\|+\|u\|_q$,
which is then equivalent to the $\ell^2$-norm $\|\cdot\|$ by \eqref{2.1}.
Therefore, we actually have $V=H=V^*$.
In accord with \eqref{2.4n}, we let
\be\label{3.59z2}\tilde{F}(t,u)=-A^2u-Af(u)-\lam u+g(t,u)\hs\mb{and}
\hs\tilde{G}(t,u)=ku+h(t).\ee

For the following discussion, we set one more hypothesis for $g(t,u)$ as follows,

\noindent\textbf{(G3)} there is $\psi_5(t)=(\psi_{5i}(t))_{i\in\Z}$, such that
\be\label{3.58z}\left|\frac{\pa g_i}{\pa u_i}(t,u_i)\right|
\leqslant\psi_{5i}(t)\mb{ with }\psi_5\in L_{\rm loc}^\infty(\R;\ell^\infty).\ee

We take $\cT$ as the class of test functions presented in Subsection \ref{ss2.5}.
To prove the final conclusion, we need the following results.

\bl\label{le3.11} Let the hypotheses \textbf{(F)}, \textbf{(G1)}, \textbf{(G2)},
\textbf{(H1)} and \textbf{(H2)} hold.
Let $\cA$ be the pullback random attractor of $\vp$
with respect to $\cD$ given in Theorem \ref{th3.6}.
Then for each $\tau\in\R$, $\bP$-a.s. $\w\in\W$, and $s,t\in\R$ with $s\leqslant t$,
the set
\be\label{5.1a}\cU[s,t]:=\Cup_{\sig\in[s,t]}\cA(\sig,\theta_{\sig-\tau}\w)\ee
is compact in $\ell^2$.\el
\bo We pick a sequence $\{\phi_k\}_{k\in\N^+}$ in $\cU[s,t]$
arbitrarily with $\phi_k\in\cA(\sig_k,\theta_{\sig_k-\tau}\w)$,
and it suffices to show that this sequence has a convergent subsequence
with its limit in $\cU[s,t]$.

Since $\sig_k\in[s,t]$, we can assume that $\sig_k\ra\sig_0\in[s,t]$ as $k\ra\8$
(selecting its subsequence if necessary).
Take
$$\sig_*=\inf_{k\in\N}\sig_k\in[s,t].$$
Then by the invariance of the attractor, we have
$\tilde{\phi}_k\in\cA(\sig_*,\theta_{\sig_*-\tau}\w)$ for each $n\in\N$ such that
$$\phi_k=u(\sig_k,\sig_*,\theta_{-\tau}\w,\tilde{\phi}_k).$$
Moreover, since $\tilde{\phi}_{k}\in\cA(\sig_*,\theta_{\sig_*-\tau}\w)$,
which is compact in $\ell^2$,
we can furthermore assume that
$$\tilde{\phi}_k\ra\phi_*\in\cA(\sig_*,\theta_{\sig_*-\tau}\w),\hs
\mb{as }n\ra\8.$$
By (1) of Theorem \ref{th3.3},
we can infer from $\sig_k\ra\sig_0$ and $\tilde{\phi}_k\ra\phi_*$ that
$$
\|\phi_k-u(\sig_0,\sig_*,\theta_{-\tau}\w,\phi_*)\|\leqslant
\|u(\sig_k,\sig_*,\theta_{-\tau}\w,\tilde{\phi}_k)
-u(\sig_0,\sig_*,\theta_{-\tau}\w,\phi_*)\|\ra0,$$
as $n\ra\8$.
The compactness is thus proved.
\eo

\bl\label{le3.13} Let the hypotheses \textbf{(F)}, \textbf{(G1)}, \textbf{(G2)}, \textbf{(G3)},
\textbf{(H1)} and \textbf{(H2)} hold. Let $\Psi\in\cT$ and $u\in\cA(t,\theta_{t-\tau}\w)$,
where $\cA=\{\cA(\tau,w)\}_{(\tau,\w)\in\R\X\W}$ is the pullback random attractor given
by Theorem \ref{th3.6}.
Then the mapping
$$\cA(t,\theta_{t-\tau}\w)\ni u\mapsto\langle\tilde{F}(t,u),\Psi'(u)\rangle$$
is continuous in norm of $H$, i.e., for a sequence $\{u^{(n)}\}_{n\in\N^+}$ and $u$ in $V$,
\be\label{3.59z}\langle\tilde{F}(t,u^{(n)}),\Psi'(u^{(n)})\rangle\ra
\langle\tilde{F}(t,u),\Psi'(u)\rangle,\hs\mb{as }\|u^{(n)}-u\|\ra0.\ee
And for the invariant sample measures $\{\mu_{\tau,\w}\}_{(\tau,\w)\in\R\X\W}$ given
in Theorem \ref{th3.11}, for fixed $\tau\in\R$ and $\bP$-a.s. $\w\in\W$, the mapping
\be\label{3.60z}
\vsig\mapsto\int_{H}\langle\tilde{F}(\vsig,u),\Psi'(u)\rangle
\mu_{\vsig,\theta_{\vsig-\tau}\w}(\di u)
\ee
makes sense and belongs to $L^1_{\rm loc}(\R)$.
\el

\bo We first prove the continuity of the mapping
\be\label{3.61z}V\ni u\mapsto\tilde{F}(t,u)\in V^*\ee
in norm of $H$ for each $t\in\R$.
Take $u$, $v\in V$ with $\|u\|$, $\|v\|\leqslant R$ for some $R>0$
and set $\tilde{w}=u-v$.
For each fixed $t\in\R$, we observe by \eqref{3.1z}, \eqref{3.58z} and
the Differential Mean-Value Theorem that
\begin{align*}&\|\tilde{F}(t,u)-\tilde{F}(t,v)\|\\
\leqslant&(16+\lam)\|\tilde{w}\|+4\|f(u)-f(v)\|+\|g(t,u)-g(t,v)\|\\
\leqslant&\(16+\lam\)\|\tilde{w}\|+\(\sum_{i\in\Z}\gam^2(R)|\tilde{w}_i|^2\)^{\frac12}
+\(\sum_{i\in\Z}\left|\frac{\pa g_i}{\pa u_i}(t,\xi_i)\right|^2
|\tilde{w}_i|^2\)^{\frac12}\\
\leqslant&\(16+\lam+4\gam(R)+\|\psi_5(t)\|_\8\)\|\tilde{w}\|,\end{align*}
which indicates the continuity \eqref{3.61z} in $H$ (also $V$ of course) immediately.

The condition (2) for the class $\cT$ also implies
the continuity of $u\mapsto\Psi'(u)$ from $V$ to $V$ in norm of $H$.
As a result, the continuity \eqref{3.59z} holds true.
\vs

For \eqref{3.60z}, due to the fact that for each $(\tau,\w)\in\R\X\W$,
$\cA(\tau,\w)$ is compact in $H$, and the invariant sample measures $\mu_{\tau,\w}$ is
supported by $\cA(\tau,\w)$,
the mapping \eqref{3.60z} surely makes sense.
Note also that
\begin{align}|\langle \tilde{F}(\vsig,u),\Psi'(u)\rangle|
\leqslant&|\( A^2u,\Psi'(u)\)|
+\lam\left|\( u,\Psi'(u)\)\right|+\left|\langle Af(u),\Psi'(u)\rangle\right|
+\left|\langle g(\vsig,u),\Psi'(u)\rangle\right|\notag\\
\leqslant&(16+\lam)\|u\|\|\Psi'(u)\|+4\|f(u)\|_1\|\Psi'(u)\|
+\|g(\vsig,u)\|\|\Psi'(u)\|\label{5.2a}\\
\leqslant&c\|\Psi'(u)\|_V\(\|u\|+\|\al_1\|_1\|u\|^{q-1-\epsilon}+\|\al_2\|_1
+\|\psi_3(\vsig)\|_\8\|u\|^{q-1}+\|\psi_4(\vsig)\|\).\notag
\end{align}
We consider an arbitrary interval $[s,t]\subset\R$.
By (2) of the definition of $\cT$, we know $\|\Psi'(u)\|_V$ is bounded;
by the compactness of $\cU[s,t]$, we also know that $\|u\|$ is bounded.
Hence we have a positive constant $\sM_1$ independent of $\vsig$ such that
$$|\langle \tilde{F}(\vsig,u),\Psi'(u)\rangle|
\leqslant\sM_1(1+\|\al_1\|_1+\|\al_2\|_1+\|\psi_3(\vsig)\|_\8+\|\psi_4(\vsig)\|^2)$$
and then
\begin{align*}&\left|\int_{H}\langle \tilde{F}(\vsig,u),\Psi'(u)\rangle
\mu_{\vsig,\theta_{\vsig-t}\w}(\di u)\right|
=\left|\int_{\cA(\vsig,\theta_{\vsig-t}\w)}
\langle \tilde{F}(\vsig,u),\Psi'(u)\rangle
\mu_{\vsig,\theta_{\vsig-t}\w}(\di u)\right|\\
\leqslant&\sM_1(1+\|\al_1\|_1+\|\al_2\|_1+\|\psi_3(\vsig)\|_\8+\|\psi_4(\vsig)\|^2),
\end{align*}
which ensures that \eqref{3.60z} belongs to $L^1_{\rm loc}(\R)$ by
the hypothesis \textbf{(G1)}.
The proof is complete.\eo

\bl\label{le3.14} Under the condition of Lemma \ref{le3.13}, the mapping
$\cA(t,\theta_{t-\tau}\w)\ni u\mapsto(\tilde{G}(t,u),\Psi'(u))$
is continuous in norm of $H$.
And for the invariant sample measures $\{\mu_{\tau,\w}\}_{(\tau,\w)\in\R\X\W}$ given
in Theorem \ref{th3.11}, for fixed $\tau\in\R$ and $\bP$-a.s. $\w\in\W$, the mapping
\be\label{3.60zz}
\vsig\mapsto\int_{H}(\tilde{G}(\vsig,u),\Psi'(u))\mu_{\vsig,\theta_{\vsig-\tau}\w}(\di u)
\ee
makes sense and belongs to $L^2_{\rm loc}(\R)$, which guarantees the stochastic integral
\be\label{3.60az}
\int_s^t\int_{H}(\tilde{G}(\vsig,u),\Psi'(u))
\mu_{\vsig,\theta_{\vsig-\tau}\w}(\di u)\di\vsig
\ee
$t\geqslant s$, is well defined.
\el

\bo By the discussion in the proof of Lemma \ref{le3.13},
the mapping $H\ni u\mapsto\Psi'(u)\in H$ is continuous.
And $\tilde{G}(t,u)$ is obviously continuous in $u$ in $H$, which assures the continuity of
$\cA(t,\theta_{t-\tau}\w)\ni u\mapsto(\tilde{G}(t,u),\Psi'(u))$.

For the second conclusion, pick an arbitrary interval $[s,t]\subset\R$.
We recall the condition (2) of Definition \ref{de1.1} and find that
$\Psi'(u)$ is bounded on $V$ and hence also on $H$.
Combining the compactness of $\cU[s,t]$ by Lemma \ref{le3.11}, which implies $\|u\|$ is bounded
on $\cU[s,t]$, we have a positive constant $\sM_2$ independent of $\vsig$ such that
\begin{align*}&\left|\int_{H}(\tilde{G}(\vsig,u),\Psi'(u))
\mu_{\vsig,\theta_{\vsig-\tau}\w}(\di u)\right|^2
\leqslant\sM_2\left|\int_{\cA(\vsig,\theta_{\vsig-\tau}\w)}\(1+\|h(\vsig)\|\)
\mu_{\vsig,\theta_{\vsig-\tau}\w}(\di u)\right|^2\\
=&\sM_2(1+\|h(\vsig)\|)^2\leqslant2\sM_2(1+\|h(\vsig)\|^2).
\end{align*}
Note by \textbf{(H1)} that $h\in L^2_{\rm loc}(\R,\ell^2)$.
We can then deduce that the mapping \eqref{3.60zz} belongs to $L^2_{\rm loc}(\R)$.

At last, we can refer to \cite[Section 2.3]{K05} that once \eqref{3.60zz} belongs to
$L^2_{\rm loc}(\R)$, the stochastic (Wiener) integral \eqref{3.60az} is well defined.
The proof ends here.\eo

\bl\label{le3.15} Under the condition of Lemma \ref{le3.13}, the mapping
\be\label{3.60bz}\cA(t,\theta_{t-\tau}\w)\ni
u\mapsto\Psi''(u)(\tilde{G}(t,u),\tilde{G}(t,u))\ee
is continuous in norm of $H$.
And for the invariant sample measures $\{\mu_{\tau,\w}\}_{(\tau,\w)\in\R\X\W}$ given
in Theorem \ref{th3.11}, for fixed $\tau\in\R$ and $\bP$-a.s. $\w\in\W$, the mapping
\be\label{3.60zzz}
\vsig\mapsto\int_{H}\Psi''(u)(\tilde{G}(t,u),\tilde{G}(t,u))
\mu_{\vsig,\theta_{\vsig-\tau}\w}(\di u)
\ee
makes sense and belongs to $L^1_{\rm loc}(\R)$.
\el

\bo By (3) of Definition \ref{de1.1}, we know $\Psi''(u)$ is continuous
from $H$ to $L(H\X H,\R)$.
And $\tilde{G}(t,u)=ku+h(t)$ is evidently continuous from $H$ to $H$.
Hence the continuity of the mapping \eqref{3.60bz} follows trivially.

Recall (3) of Definition \ref{de1.1} again, we know $\Psi''(u)$ is bounded
from $H$ to $L(H\X H,\R)$.
Let $s,t\in\R$ with $s\leqslant t$.
Using again the boundedness of $\|u\|$ on $\cU[s,t]$, we have $\sM_3>0$
independent of $\vsig$ such that
\begin{align*}&\left|\int_{H}\Psi''(u)(\tilde{G}(\vsig,u),\tilde{G}(\vsig,u))
\mu_{\vsig,\theta_{\vsig-\tau}\w}(\di u)\right|
\leqslant\sM_3\left|\int_{\cA(\vsig,\theta_{\vsig-\tau}\w)}\(1+\|h(\vsig)\|\)^2
\mu_{\vsig,\theta_{\vsig-\tau}\w}(\di u)\right|\\
=&\sM_3(1+\|h(\vsig)\|)^2\leqslant2\sM_3(1+\|h(\vsig)\|^2),
\end{align*}
which is in $L^1_{\rm loc}(\R)$.
We complete the proof now.
\eo

Now we are well prepared to show the final consequence for the case when $k>0$.
\bt\label{th3.15} Suppose that the hypotheses \textbf{(F)}, \textbf{(G1)}, \textbf{(G2)},
\textbf{(G3)}, \textbf{(H1)} and \textbf{(H2)} hold.
Let $\cA$ be the pullback random $\cD$-attractor of $\vp$ given in Theorem \ref{th3.6} and
$\{\mu_{\tau,\w}\}_{(\tau,\w)\in\R\X\W}$ be the invariant sample measures constructed in
Theorem \ref{th3.11}.
Then the invariant sample measures $\{\mu_{\tau,\w}\}_{(\tau,\w)\in\R\X\W}$ satisfy
the stochastic Liouville type equation \eqref{1.7} for each $\Psi\in\cT$.\et
\bo Based on the solution \eqref{3.6z} and definition of the class $\cT$,
and by \eqref{3.55z}, \eqref{3.57z} and \eqref{1.6}, we know that for all $t\geqslant\tau$,
\begin{align}&\int_{H}\Psi(\phi)\mu_{t,\theta_{t-\tau}\w}(\di\phi)
-\int_{H}\Psi(\phi)\mu_{s,\theta_{s-\tau}\w}(\di\phi)\notag\\
=&\int_{\cA(t,\theta_{t-\tau}\w)}\Psi(\phi)\mu_{t,\theta_{t-\tau}\w}(\di\phi)
-\int_{\cA(s,\theta_{s-\tau}\w)}\Psi(\phi)\mu_{s,\theta_{s-\tau}\w}(\di\phi)\notag\\
=&\int_{\cA(s,\theta_{s-\tau}\w)}\(\Psi(u(t,s,\theta_{s-\tau}\w,\phi))-\Psi(\phi)\)
\mu_{s,\theta_{s-\tau}\w}(\di\phi)\notag\\
=&\int_{\cA(s,\theta_{s-\tau}\w)}\int_s^t\langle\tilde{F}(\vsig,u(\vsig)),
\Psi'(u(\vsig))\rangle\di\vsig\cdot\mu_{s,\theta_{s-\tau}\w}(\di\phi)\notag\\
&+\int_{\cA(s,\theta_{s-\tau}\w)}\int_{s}^t
\(\tilde{G}(\vsig,u(\vsig)),\Psi'(u(\vsig))\)\di\tilde{W}(\vsig)\cdot
\mu_{s,\theta_{s-\tau}\w}(\di\phi)\notag\\
&+\frac12\int_{\cA(s,\theta_{s-\tau}\w)}\int_{s}^t\Psi''(u(\vsig))
\(\tilde{G}(\vsig,u(\vsig)),\tilde{G}(\vsig,u(\vsig))\)\di\vsig\cdot
\mu_{s,\theta_{s-\tau}\w}(\di\phi)\notag\\
:=&M_1+M_2+M_3,\label{3.62z}
\end{align}
where $\tilde{W}(\vsig)=W(-\tau+\vsig)-W(-\tau)$ comes from
the stochastic datum $\theta_{-\tau}\w$.
Note by (3) of Definition \ref{de2.1} and \eqref{3.24z} that
\be\label{3.63z}
u(\vsig,s,\theta_{-\tau}\w,u(s,\sig,\theta_{-\tau}\w,\phi))=
u(\vsig,\sig,\theta_{-\tau}\w,\phi).\ee
Then using \eqref{3.55z}, \eqref{3.63z}, Fubini's Theorem and
the invariance of $\mu_{\tau,\w}$ to $M_1$, we have
\begin{align}M_1=&\LIM_{\tau'\ra-\8}\frac{1}{s-\tau'}\int_{\tau'}^{s}\int_{H}\int_s^t
\langle\tilde{F}(\vsig,u(\vsig,\sig,\theta_{-\tau}\w,\phi)),
\Psi'(u(\vsig,\sig,\theta_{-\tau}\w,\phi))\rangle
\di\vsig\cdot \mu_{\sig,\theta_{\sig-\tau}\w}(\di\phi)\di\sig\notag\\
=&\LIM_{\tau'\ra-\8}\frac{1}{s-\tau'}\int_{\tau'}^{s}\int_s^t\int_{H}
\langle\tilde{F}(\vsig,u(\vsig,\sig,\theta_{-\tau}\w,\phi)),
\Psi'(u(\vsig,\sig,\theta_{-\tau}\w,\phi))\rangle
\mu_{\sig,\theta_{\sig-\tau}\w}(\di\phi)\di\vsig\di\sig\notag\\
=&\LIM_{\tau'\ra-\8}\frac{1}{s-\tau'}\int_{\tau'}^{s}\int_s^t\int_{H}
\langle\tilde{F}(\vsig,u),\Psi'(u)\rangle
\mu_{\vsig,\theta_{\vsig-\tau}\w}(\di u)\di\vsig\di\sig\notag\\
=&\int_s^t\int_{H}\langle\tilde{F}(\vsig,u),\Psi'(u)\rangle
\mu_{\vsig,\theta_{\vsig-\tau}\w}(\di u)\di\vsig,\label{3.64z}
\end{align}
where we have replaced $t$, $\tau$ and $\w$ in \eqref{3.57z} by $\vsig-\sig$, $\sig$
and $\theta_{\sig-\tau}\w$, respectively, to obtain the third equality.
Similarly, for $M_2$ and $M_3$, we also have
\begin{align}M_2=&\LIM_{\tau'\ra-\8}\frac{1}{s-\tau'}\suo\int_{\tau'}^{s}\suo\int_{H}\!\!\int_s^t
\suo\(\tilde{G}(\vsig,u(\vsig,\sig,\theta_{-\tau}\w,\phi)),
\Psi'(u(\vsig,\sig,\theta_{-\tau}\w,\phi))\)
\di\tilde{W}(\vsig)\cdot \mu_{\sig,\theta_{\sig-\tau}\w}(\di\phi)\di\sig\notag\\
=&\int_s^t\int_{H}\(\tilde{G}(\vsig,u),\Psi'(u)\)
\mu_{\vsig,\theta_{\vsig-\tau}\w}(\di u)\di\tilde{W}(\vsig)\hs\mb{and}
\end{align}
\be M_3=\frac12\int_s^t\int_{H}\Psi''(u)\(\tilde{G}(\vsig,u),\tilde{G}(\vsig,u)\)
\mu_{\vsig,\theta_{\vsig-\tau}\w}(\di u)\di \vsig.\label{3.66z}
\ee
Then the stochastic Liouville type equation \eqref{1.7} can be obtained by combining
the equalities from \eqref{3.62z} to \eqref{3.66z}.
The proof is thus finished.
\eo

According to the discussion above, we naturally obtain
the existence of sample statistical solutions.
\bt Suppose that the hypotheses \textbf{(F)}, \textbf{(G1)}, \textbf{(G2)},
\textbf{(G3)}, \textbf{(H1)} and \textbf{(H2)} hold.
Then the family of Borel probability measures
$\{\mu_{\tau,\w}\}_{(\tau,\w)\in\R\X\W}$ obtained
in Theorem \ref{th3.11} is a sample statistical solution
of the equation \eqref{2.4n}.
\et
\bo We check the conditions of Definition \ref{de2.9n} one by one.
For the condition (1), we let $\Psi\in\cC_b(H)$.
Then by \eqref{3.55z} and \eqref{3.57z}, for $t,\,s,\,s_*\in\R$ with $t,\,s>s_*$
and $s,\,s_*$ fixed, we have
\begin{align}&\left|\int_{H}\Psi(u)\mu_{t,\theta_{t-\tau}\w}(\di u)
-\int_H\Psi(u)\mu_{s,\theta_{s-\tau}\w}(\di u)\right|\notag\\
\leqslant&\int_{\cA(s_*,\theta_{s_*-\tau}\w)}
\left|\Psi(u(t,s_*,\theta_{-\tau}\w,\phi))-\Psi(u(s,s_*,\theta_{-\tau}\w,\phi))\right|
\mu_{s_*,\theta_{s_*-\tau}\w}(\di\phi).\label{3.70z}
\end{align}
Since $\cA(s_*,\theta_{s_*-\tau}\w)$ is compact,
by (1) of Theorem \ref{th3.3}, we know the continuity of
$s\mapsto u(s,s_*,\theta_{-\tau}\w,\phi)$ is uniform for all $\phi\in\cA(s_*,\theta_{s_*-\tau}\w)$.
Hence for every $\de>0$, there exists $\ol{\de}>0$
independent of $\phi\in\cA(s_*,\theta_{s_*-\tau}\w)$, such that when $|t-s|<\ol{\de}$,
\be\label{3.71z}\|u(t,s_*,\theta_{-\tau}\w,\phi)-u(s,s_*,\theta_{-\tau}\w,\phi)\|<\de,\ee
for all $\phi\in\cA(s_*,\theta_{s_*-\tau}\w)$.
Furthermore, for all $\phi\in\cA(s_*,\theta_{s_*-\tau}\w)$ and $\vsig\in[s-1,s+1]$,
we know that
$$u(\vsig,s_*,\theta_{-\tau}\w,\phi)\in\cU[s-1,s+1],$$
which is compact by Lemma \ref{le3.11}.
Also, $\Psi\in\cC_b(H)$ implies that $\Psi$ is uniformly continuous on $\cU[s-1,s+1]$.
This means that, for each $\ve>0$, there is $\de>0$ such that when $u,\,v\in\cU[s,t]$
and $\|u-v\|<\de$,
\be\label{3.72z}|\Psi(u)-\Psi(v)|<\ve.\ee

As a result, by \eqref{3.70z}, \eqref{3.71z} and \eqref{3.72z}, it can be seen that
for each $\ve>0$, there exists $\de>0$ and $\ol{\de}\in(0,1)$, such that when $|t-s|<\ol{\de}$,
\eqref{3.71z} holds true and then
$$\eqref{3.70z}<\int_{\cA(s_*,\theta_{s_*-\tau}\w)}\ve
\mu_{s_*,\theta_{s_*-\tau}\w}(\di\phi)=\ve,$$
which indicates the condition (1) of Definition \ref{de2.9n}.

For (2) of Definition \ref{de2.9n},
we can recall Lemmas \ref{le3.13}, \ref{le3.14} and \ref{le3.15}.
The conclusions for $\langle\tilde{F}(s,u),\phi\rangle$ and $(\tilde{G}(s,u),\phi)$
can be obtained by defining $\Psi(u)\equiv\phi$ in Lemmas \ref{le3.13} and \ref{le3.14}.
The conclusions for $\Phi(\tilde{G}(s,u),\tilde{G}(s,u))$ can be obtained by
defining $\Psi(u)=\Phi(u,u)$.

The condition (3) of Definition \ref{de2.9n} follows immediately from Theorem \ref{th3.15}.
The theorem has been thus proved here.
\eo
\section{The case when $k=0$}

In this section, we discuss the problem \eqref{2.4n} for the case when $k=0$.
Still, we impose the hypotheses \textbf{(F)}, \textbf{(G1)}, \textbf{(G2)}, \textbf{(G3)},
\textbf{(H1)} and \textbf{(H2)} where necessary on $f$, $g$ and $h$ in \eqref{2.4n}.

\subsection{The induced nonautonomous random dynamical system}

We still use the concepts of local and global solutions given in Subsection \ref{ss3.1}.
In order to obtain the unique existence of solutions, we let $u(t,\tau,\w,u_\tau)$ be
a solution of \eqref{2.4n} and set
$$v(t)=u(t)-z(\theta_t\w)h(t),\Hs
v(\tau)=v_\tau:=u_\tau-z(\theta_\tau\w)h(\tau),$$
where $z(\theta_t\w)$ is given in \eqref{2.5n}.
Then similar to \eqref{3.8z}, we can deduce that $v(t)$ satisfies
\be\label{4.1}\frac{\di v(t)}{\di t}=\ol{F}(t,\w,v),\hs t>\tau,\Hs v(\tau)=v_\tau,\ee
where
\begin{align*}\ol{F}(t,\w,v):=&(-A^2 v-\lam v)
+\(-Af(v+z(\theta_t\w)h(t))+g(t,v+z(\theta_t\w)h(t))\)\\
&+z(\theta_t\w)\(-A^2h(t)+(1-\lam)h(t)-h'(t)\)\\
:=&\ol{F}_1(v)+\ol{F}_2(t,\w,v)+\ol{F}_3(t,\w).
\end{align*}

Now we prove the local unique existence of solutions for \eqref{4.1}
in the following theorem.

\bt Let the hypotheses \textbf{(F)}, \textbf{(G1)} and \textbf{(H1)} hold.
Then for each $(\tau,v_\tau)\in\R\X\ell^2$ and $\bP$-a.s. $\w\in\W$,
the system \eqref{4.1} has a unique global solution
$$v(\cdot,\tau,\w,v_\tau)\in\cC^1([\tau,\8);\ell^2)\cap L^q_{\rm loc}(0,\8;\ell^q)$$
such that $v(\tau,\tau,\w,v_\tau)=v_\tau$.
\et
\bo Like before, we first consider the local solutions.
Let $T>\tau$ and $t\in[\tau,T]$ and pick $u,v\in\ell^2$ such that $\|u\|$, $\|v\|<R$.
Notice the following estimates:
\be\label{4.2}
\(\ol{F}_1(u)-\ol{F}_1(v),u-v\)=-\|A(u-v)\|^2-\lam\|u-v\|^2\hs\mb{and}
\ee
\begin{align}
\(\ol{F}_2(t,\w,u)-\ol{F}_2(t,\w,v),u-v\)\leqslant&\left|\(A\(f(u+z(\theta_t\w)h(t))
-f(v+z(\theta_t\w)h(t))\),u-v\)\right|\notag\\&
+\(g(t,u+z(\theta_t\w)h(t))-g(t,v+z(\theta_t\w)h(t)),u-v\)\notag\\
\leqslant&(4\gamma(R)+\|\psi_2(t)\|_\8)\|u-v\|^2.\label{4.3}
\end{align}
We can then deduce by \eqref{4.2} and \eqref{4.3} that
\begin{align}\(\ol{F}(t,\w,u)-\ol{F}(t,\w,v),u-v\)\leqslant
&\(\esssup_{t\in[\tau,T]}\|\psi_2(t)\|_\8+4\gamma(R)
\)\|u-v\|^2\notag\\
:=&\ol{\cP}_1\(\tau,T,R\)\|u-v\|^2,\label{4.4}\end{align}
which ensures the local unique existence of solutions for \eqref{4.1}.

Next we prove the global existence of solutions for the problem \eqref{4.1}.
Observe similarly that
\begin{align}\label{4.5}(\ol{F}_1(t,\w,v),v)=&-\|Av\|^2-\lam\|v\|^2,\\
(\ol{F}_2(t,\w,v),v)=&\(\ol{F}_2(t,\w,v),v+z(\theta_t\w)h(t)\)
-\(\ol{F}_2(t,\w,v),z(\theta_t\w)h(t)\)\notag\\
:=&\ol{I}_1+\ol{I}_2,\end{align}
\begin{align}\ol{I}_1\leqslant&
\sum_{i\in\Z}\(\al_{1i}|v_i+z(\theta_t\w)h_i(t)|^{q-1-\epsilon}+\al_{2i}\)
\left|\(A(v+z(\theta_t\w)h(t)\)_i\right|\notag\\
&-\beta\|v+z(\theta_t\w)h(t)\|_q^q+\|\psi_1(t)\|_1\notag\\
\leqslant&c\(\|\al_1\|_1^{\frac{q}{\epsilon}}+\|\al_2\|_1^{\frac{q}{q-1}}\)
+\|\psi_1(t)\|_1-\frac{3\beta}{4}\|v+z(\theta_t\w)h(t)\|_q^q,
\end{align}
\begin{align}|\ol{I}_2|\leqslant&
\left|\(f(v+z(\theta_t\w)h(t)),z(\theta_t\w)Ah(t)\)\right|
+\left|\(g(t,v+z(\theta_t\w)h(t),z(\theta_t\w)h(t)\)\right|\notag\\
\leqslant&4\(\|\al_1\|_1\|v+z(\theta_t\w)h(t)\|_q^{q-1-\epsilon}
+\|\al_2\|_1\)\|z(\theta_t\w)h(t)\|_q\notag\\
&+\sum_{i\in\Z}|\psi_{3i}(t)||v+z(\theta_t\w)h_i(t)|^{q-1}|z(\theta_t\w)h_i(t)|
+\|\psi_4(t)\|\|z(\theta_t\w)h(t)\|\notag\\
\leqslant&c\(\|\al_1\|_1^{\frac{q}{\epsilon}}+\|\al_2\|_1^{\frac{q}{q-1}}+\|\psi_4(t)\|^2\)
+c\(1+\|\psi_3(t)\|^q_\8\)|z(\theta_t\w)|^q\|h(t)\|_q^q\notag\\
&+c|z(\theta_t\w)|^2\|h(t)\|^2+\frac{\beta}{4}\|v+z(\theta_t\w)h(t)\|_q^q\hs\mb{and}
\end{align}
\be\label{4.9}|(\ol{F}_3(t,\w),v)|\leqslant
c|z(\theta_t\w)|^2\(\|h(t)\|^2+\|h'(t)\|^2\)+\frac{\lam}{2}\|v\|^2.
\ee
Then considering the inner product of \eqref{4.1} and $v$,
we obtain by the estimates from \eqref{4.5} and \eqref{4.9} that
\be\frac{\di\|v(t)\|^2}{\di t}+\lam\|v\|^2
+\beta\|v+z(\theta_t\w)h(t)\|_q^q
\leqslant c\ol{\cQ}_1(t,\w),\label{4.10}
\ee
where
\begin{align*}\ol{\cQ}_1(t,\w):=&\(\|\al_1\|_1^{\frac{q}{\epsilon}}+\|\al_2\|_1^{\frac{q}{q-1}}
+\|\psi_1(t)\|_1+\|\psi_4(t)\|^2\)+|z(\theta_t\w)|^2\(\|h(t)\|^2+\|h'(t)\|^2\)\\
&+\(1+\|\psi_3(t)\|_\8^q\)|z(\theta_t\w)|^q\|h(t)\|_q^q.
\end{align*}

Utilizing Gronwall's lemma to \eqref{4.10} over the interval $[\sig,s]\subset[\tau,T)$,
we obtain
\be\|v(s)\|^2+\beta\int_\sig^s\me^{\lam(\vsig-s)}
\|v(\vsig)+z(\theta_\vsig\w)h(\vsig)\|_q^q\di\vsig
\leqslant\me^{\lam(\sig-s)}\|v(\sig)\|^2
+c\int_{\sig}^s\me^{\lam(\vsig-s)}\ol{\cQ}_1(\vsig,\w)\di\vsig.\label{4.11}
\ee
Replacing $\sig$ and $s$ by $\tau$ and $t$ respectively for $t\in[\tau,T)$
and by the hypotheses \textbf{(F)}, \textbf{(G1)} and \textbf{(H1)}, we have
\begin{align}&\|v(t)\|^2+\beta\int_\tau^t\me^{\lam(\vsig-t)
}\|v(\vsig)+z(\theta_\vsig\w)h(\vsig)\|_q^q\di\vsig\notag\\
\leqslant&\|v_\tau\|^2
+c\int_{\tau}^T\ol{\cQ}_1(\vsig,\w)\di\vsig
:=\ol{\cP}_2:=\ol{\cP}_2(\tau,T,\w,\|v_\tau\|)<\8.\label{4.12}
\end{align}

In order to extend the local solution to a global one,
we still only need to show that $v(t_n)$ is a Cauchy sequence in $\ell^2$
for an arbitrary strictly increasing sequence $\{t_n\}_{n\in\N^+}$ in $(\tau,T)$
with $t_n\ra T^-$.
With the same analysis, we can analogously obtain
$$\|v(t_n)-v(t_m)\|\leqslant\|v(\tau+t_n-t_m)-v(\tau)\|
\me^{2\ol{\cP}_1\(\tau,T,\sqrt{\ol{\cP}_2}\)(T-\tau)}.$$
Hence $v(t_n)$ is a Cauchy sequence in $\ell^2$.

Finally, by the inequality \eqref{4.11} and \eqref{4.12}, we know that
\begin{align}\int_{\tau}^T\|v(s)\|_q^q\di s\leqslant&
2^q\(\int_{\tau}^T\|v(s)+z(\theta_s\w)h(s)\|_q^q\di s
+\int_\tau^T\|z(\theta_s\w)h(s)\|_q^q\di s\)\notag\\
\leqslant&{2^q}\me^{\lam(T-\tau)}
\int_{\tau}^{T}\me^{\lam(s-T)}\|v(s)+z(\theta_s\w)h(s)\|_q^q\di s\notag\\
&+2^q\int_\tau^T\|z(\theta_s\w)h(s)\|_q^q\di s<\8,\end{align}
which means that $v(\cdot,\tau,\w,v_\tau)\in L^q_{\rm loc}(0,\8;\ell^q)$.
The proof is complete.
\eo

Then we set
\be\label{4.14}
u(t,\tau,\w,u_\tau)=v(t,\tau,\w,u_\tau-z(\theta_\tau\w)h(\tau))+z(\theta_t\w)h(t).
\ee
The mapping $u(t)$ in \eqref{4.14} thus defines the unique global solution of \eqref{4.1}
for each $(\tau,u_\tau)\in\R\X\ell^2$ and $\bP$-a.s. $\w\in\W$.
\Vs

Similar to the treatment in Subsection \ref{ss3.2},
we define a mapping $\vp:\R^+\X\R\X\W\X\ell^2\ra\ell^2$,
such that for every $(t,\tau,\w,u_\tau)\in\R^+\X\R\X\W\X\ell^2$,
\be\label{4.15}\vp(t,\tau,\w,u_\tau)=u(t+\tau,\tau,\theta_{-\tau}\w,u_\tau)
=v(t+\tau,\tau,\theta_{-\tau}\w,v_\tau)+z(\theta_{t+\tau}\w)h(t+\tau)\ee
with $v_\tau:=u_\tau-z(\theta_\tau\w)h(\tau)$.
With almost the same proof of Theorem \ref{th3.3}, we know that
$\vp$ defined in \eqref{4.15} induces a nonautonomous random dynamical system.

\subsection{Existence of pullback random attractors}

In this subsection, we consider the existence of pullback random attractors for
the NRDS $\vp$ in \eqref{4.15}.
Here we still pick $\lam_0$, $\lam_1\in\R$ with $0<\lam_1<\lam_0<\lam$.
However, in the new situation, we take the nonautonomous random set
$D=\{D(\tau,\w):\tau\in\R,\w\in\W\}$ in $\ell^2$ such that
for every $\tau\in\R$ and $\w\in\W$,
\be\label{4.16}\lim_{s\ra -\8}\me^{\lam_0s}\|D(\tau+s,\theta_{s}\w)\|^2=0,\ee
and let $\ol{\cD}$ be the universe of all the nonautonomous random sets $D$ satisfying
\eqref{4.15}.

\bl Let the hypotheses \textbf{(F)}, \textbf{(G1)}, \textbf{(G2)},
\textbf{(H1)} and \textbf{(H2)} hold and
$v_{\tau-t}+z(\theta_{-t}\w)h(\tau-t)\in D(\tau-t,\theta_{-t}\w)$ with $D\in\ol{\cD}$.
Then for $\bP$-a.s. $\w\in\W$, there exists a certain $\ol{T}_0=\ol{T}_0(\tau,\w,D)>0$ such that
when $t>\ol{T}_0$, it holds that
\be\label{4.17}\|v(\tau,\tau-t,\theta_{-\tau}\w,v_{\tau-t})\|^2
+\beta\suo\int_{-t}^0\me^{\lam s}\|
v(s+\tau)+z(\theta_{s}\w)h(s+\tau)\|_q^q\di s\leqslant\ol{\cR}^2(\lam,\tau,\w),\ee
where $\ol{\cR}(\lam',\tau,\w)$ is defined for each $\lam'\in[\lam_0,\lam]$ as
\be\label{4.18}
\ol{\cR}^2(\lam',\tau,\w):=1+c\int_{-\8}^0\me^{\lam' s}\ol{\cQ}_1(s+\tau,\theta_{-\tau}\w)\di s<+\8.
\ee
Moreover, for each $(\tau,\w)\in\R\X\W$, let
\be\label{4.19}
\ol{K}(\tau,\w):=\{u\in\ell^2:\|u\|\leqslant\ol{\cR}(\lam,\tau,\w)+|z(\w)|\|h(\tau)\|\}.\ee
Then
\be\label{4.20}
\lim_{s\ra-\8}\me^{\lam_0 s}\|\ol{K}(\tau+s,\theta_s\w)\|^2=0.
\ee
\el
\bo Let $\tau-t$, $\tau-t$ and $\theta_{-\tau}\w$ take place of
$\tau$, $\sig$ and $\w$ in \eqref{4.11}, respectively.
We obtain that,
\begin{align}&\|v(s,\tau-t,\theta_{-\tau}\w,v_{\tau-t})\|^2
+\beta\int_{-t}^{s-\tau}\me^{\lam(\vsig+\tau-s)}
\|v(\vsig+\tau)+z(\theta_{\vsig}\w)h(\vsig+\tau)\|_q^q\di\vsig\notag\\
\leqslant&\me^{\lam(\tau-t-s)}\|v_{\tau-t}\|^2
+c\int_{-t}^{s-\tau}\me^{\lam(\vsig+\tau-s)}
\cQ_1(\vsig+\tau,\theta_{-\tau}\w)\di\vsig.
\label{4.21}\end{align}
If we choose $s=\tau$ and $\lam'\in[\lam_0,\lam)$, the inequality \eqref{4.21} becomes
\begin{align}&\|v(\tau,\tau-t,\theta_{-\tau}\w,v_{\tau-t})\|^2
+\beta\int_{-t}^{0}\me^{\lam s}\|v(s+\tau)+z(\theta_{s}\w)h(s+\tau)\|_q^q\di s\notag\\
\leqslant&\me^{-\lam' t}\|v_{\tau-t}\|^2
+c\int_{-t}^0\me^{\lam' s}\cQ_1(s+\tau,\theta_{-\tau}\w)\di s.\label{4.22}
\end{align}

It follows from \eqref{2.7n} that for all $n\in\N^+$,
\be\label{4.23}\lim_{t\ra+\8}\me^{(\lam_1-\lam')t}|z(\theta_{-t}\w)|^n
=\lim_{t\ra+\8}\(t^n\me^{(\lam_1-\lam')t}\)\frac{|z(\theta_{-t}\w)|^n}{t^n}=0.\ee
By the fact $v_{\tau-t}+h(\tau-t)\in D(\tau-t,\theta_{-t}\w)$ with $D\in\ol{\cD}$,
\eqref{4.23} and \eqref{3.27z}, we obtain that
\begin{align}\label{4.25m}\me^{-\lam't}\|v_{\tau-t}\|^2\leqslant&2\me^{-\lam't}
\(\|v_{\tau-t}+z(\theta_{-t}\w)h(\tau-t)\|^2+z^2(\theta_{-t}\w)\|h(\tau-t)\|^2\)\notag\\
\leqslant&2\me^{-\lam_0t}\|D(\tau-t,\theta_{-t}\w)\|^2
+2\(\me^{(\lam_1-\lam')t}z^2(\theta_{-t}\w)\)\!\!\(\me^{-\lam_1t}\|h(\tau-t)\|^2\)\ra0,
\end{align}
as $t\ra+\8$.
The conclusion \eqref{4.17} follows immediately from \eqref{4.25m} and \eqref{4.22}.
The boundedness in \eqref{4.18} is a deduction of the hypotheses \textbf{(G2)}, \textbf{(H2)}
and \eqref{4.23}.

For \eqref{4.20}, we recall \eqref{4.18}, \eqref{4.19}, \eqref{4.23} and \eqref{4.25m}, and have
\begin{align*}&\lim_{s\ra-\8}\me^{\lam_0 s}\|\ol{K}(\tau+s,\theta_s\w)\|^2\\
\leqslant&2\lim_{s\ra-\8}\me^{\lam_0s}\(\ol{\cR}^2(\lam,\tau+s,\theta_s\w)
+z^2(\theta_{s}\w)\|h(\tau+s)\|^2\)\\
\leqslant&c\lim_{s\ra-\8}\(\me^{\lam_0s}+\me^{\lam_0 s}
\int_{-\8}^0\me^{\lam\vsig}\ol{\cQ}_1(\vsig+\tau+s,\theta_{-\tau}\w)\di\vsig\)
+2\lim_{s\ra-\8}\me^{\lam_0s}z^2(\theta_{s}\w)\|h(\tau+s)\|^2\\
=&c\lim_{s\ra-\8}\int_{-\8}^s\me^{\lam_0\vsig}\ol{\cQ}_1(\vsig+\tau,\theta_{-\tau}\w)\di\vsig
+2\lim_{s\ra-\8}\me^{\lam_0s}z^2(\theta_{s}\w)\|h(\tau+s)\|^2=0,
\end{align*}
by \eqref{4.18}, \eqref{4.23} and \eqref{3.27z}.
The proof is finished.\eo

\bl\label{le4.3} Let the hypotheses \textbf{(F)}, \textbf{(G1)}, \textbf{(G2)},
\textbf{(H1)} and \textbf{(H2)} hold and
$v_{\tau-t}+z(\theta_{-t}\w)h(\tau-t)\in D(\tau-t,\theta_{-t}\w)$ with $D\in\ol{\cD}$.
Then for every $\ve>0$ and $\bP$-a.s. $\w\in\W$, there exist $\ol{T}_1=\ol{T}_1(\ve,\tau,\w,D)>0$
and $\ol{I}=\ol{I}(\ve,\tau,\w,D)\in\N$ such that when $t>\ol{T}_1$, it holds that
\be\sum_{|i|>\ol{I}}|
v_i(\tau,\tau-t,\theta_{-\tau}\w,v_{\tau-t})|^2<\ve.
\label{4.25}\ee
\el
\bo Let $N$ be a positive integer with $N\geqslant3$ and $w(t):=\rho^N\otimes v(t)$
for a solution $v(t,\tau,\w,v_\tau)$ of \eqref{3.8z} with $\rho^N$ given in \eqref{3.38y}.
We similarly consider $\|w(\tau,\tau-t,\theta_{-\tau}\w,v_{\tau-t})\|^2$ by noting that
$$\|w(\tau,\tau-t,\theta_{-\tau}\w,v_{\tau-t})\|^2\geqslant
\sum_{|i|>2N}|v_i(\tau,\tau-t,\theta_{-\tau}\w,v_{\tau-t})|^2.$$

We first estimate $(F(t,\w,v),w)$.
By \eqref{3.39y}, we infer that
$$
(\ol{F}_1(v),w)\leqslant\frac{8c_0}{N}\|v\|^2-\lam\|w\|^2.
$$
As to $\ol{F}_2(t,\w,v)$, we know that
\begin{align*}\(\ol{F}_2(t,\w,v),w\)=&
\(-\rho^N\otimes Af(v+z(\theta_t\w)h(t)),
\rho^N\otimes(v+z(\theta_t\w)h(t))\)\\
&+\(g(t,v+z(\theta_t\w)h(t)),
\rho^N\otimes(v+z(\theta_t\w)h(t))\)\\
&+(-z(\theta_t\w))\(-\rho^N\otimes Af(v+z(\theta_t\w)h(t)),
\rho^N\otimes h(t))\)\\
&+(-z(\theta_t\w))\(g(t,v+z(\theta_t\w)h(t)),
\rho^N\otimes h(t)\)\notag\\
:=&\ol{J}_1+\ol{J}_2+\ol{J}_3+\ol{J}_4.\end{align*}
By the hypotheses \textbf{(F)}, \textbf{(G1)}, \textbf{(H1)}, \eqref{3.40z},
\eqref{2.2n} and Young's inequality, we know that
\begin{align*}|\ol{J}_1|\leqslant&
\|\rho^N\otimes Af(v+z(\theta_t\w)h(t))\|_1
\|v+z(\theta_t\w)h(t)\|_q\\
\leqslant&4\|\rho^{N-1}\otimes f(v+z(\theta_t\w)h(t))\|_1
\|v+z(\theta_t\w)h(t)\|_q\\
\leqslant&c\(1+\|v+z(\theta_t\w)h(t)\|_q^q\)
\(\|\rho^{N-1}\otimes\alpha_{1}\|_1+\|\rho^{N-1}\otimes\alpha_{2}\|_1\),\\
\ol{J}_2\leqslant&-\beta\|\rho^N\otimes(v+z(\theta_t\w)h(t))\|_q^q
+\|\rho^N\otimes\psi_1(t)\|_1,\\
|\ol{J}_3|\leqslant&
c|z(\theta_t\w)|\|\rho^{N-1}\otimes f(v+z(\theta_t\w)h(t))\|_1\|\rho^N\otimes h(t)\|_q\\
\leqslant&c|z(\theta_t\w)|\|\rho^{N-1}\otimes\alpha_{1}\|_1\(1+\|h(t)\|_q^q
+\|v+z(\theta_t\w)h(t)\|_q^q\)\\
&+c|z(\theta_t\w)|\|\rho^{N-1}\otimes\alpha_2\|_1\(1+\|h(t)\|_q^q\)\hs\mb{and}\\
|\ol{J}_4|\leqslant&\beta\|\rho^N\otimes(v+z(\theta_t\w)h(t))\|_q^q
+c\|\rho^N\otimes\psi_4(t)\|^2\\
&+c|z(\theta_t\w)|^q\|\psi_3(t)\|^q_{\8}\|\rho^N\otimes h(t)\|_q^q
+c|z(\theta_t\w)|^2\|\rho^N\otimes h(t)\|^2.
\end{align*}
For $\ol{F}_3(t,\w)$, we have
\begin{align*}|(\ol{F}_3(t,\w),w)|\leqslant&c|z(\theta_t\w)|^2\(\|\rho^N\otimes A^2h(t)\|^2
+\|\rho^N\otimes h(t)\|^2+\|\rho^N\otimes h'(t)\|^2\)+\frac{\lam}{2}\|w\|^2\\
\leqslant&c|z(\theta_t\w)|^2\(\|\rho^{N-2}\otimes h(t)\|^2
+\|\rho^N\otimes h'(t)\|^2\)+\frac{\lam}{2}\|w\|^2.
\end{align*}

Let
\begin{align*}&\hs\ol{\cQ}^N(t,\w,v):=\\
&(1+|z(\theta_t\w)|)\(\|\rho^{N-1}\otimes\al_1\|_1+\|\rho^{N-1}\otimes\al_2\|_1\)
\(1+\|h(t)\|_q^q+\|v+z(\theta_t\w)h(t)\|_q^q\)\notag\\
&+\|\rho^N\otimes\psi_1(t)\|_1+\|\rho^N\otimes\psi_4(t)\|^2
+|z(\theta_t\w)|^q\|\psi_3(t)\|^q_{\8}\|\rho^N\otimes h(t)\|_q^q\\
&+|z(\theta_t\w)|^2\(\|\rho^{N-2}\otimes h(t)\|^2+\|\rho^N\otimes h'(t)\|^2\).
\end{align*}
Now taking the inner product of \eqref{4.1} and $w$ and introducing the estimates above,
we obtain
\be\label{4.26}\frac{\di\|w(t)\|^2}{\di t}+\lam\|w(t)\|^2\leqslant
\frac{16c_0}{N}\|v(t)\|^2+c\ol{\cQ}^N(t,\w,v(t)).\ee
Applying Gronwall's lemma to \eqref{4.26} over $[\sig,s]$ and
substituting $\tau$, $\sig$, $s$ and $\w$ by $\tau-t$, $\tau-t$, $\tau$
and $\theta_{-\tau}\w$, respectively, we have
\begin{align*}
&\|w(\tau,\tau-t,\theta_{-\tau}\w,v_{\tau-t})\|^2\\
\leqslant&\me^{-\lam t}\|v_{\tau-t}\|^2
+\frac{16c_0}{N}\int_{-t}^0\me^{\lam\vsig}
\|v(\vsig+\tau)\|^2\di\vsig+c\int_{-t}^0\me^{\lam\vsig}
\ol{\cQ}^N(\vsig+\tau,\theta_{-\tau}\w,v(\vsig+\tau))\di \vsig\\
:=&\ol{L}_1+\ol{L}_2+\ol{L}_3.\end{align*}
By \eqref{4.25m}, $\ol{L}_1$ vanishes as $t\ra+\8$.
Similar to $L_2$ discussed in the proof of Lemma \ref{le3.5} and by \eqref{4.21},
we know that when $t$ is sufficient large (independent of $\vsig\in[-t,0]$),
\begin{align*}\|v(\vsig+\tau,\tau-t,\theta_{-\tau}\w,v_{\tau-t})\|^2
\leqslant&\me^{-\lam(t+\vsig)}\|v_{\tau-t}\|^2+c\int_{-t}^{\vsig}\me^{\lam(\sig-\vsig)}
\ol{\cQ}_1(\sig+\tau,\theta_{-\tau}\w)\di\sig\\
\leqslant&\me^{-\lam_0\vsig}\(\me^{-\lam_0t}\|v_{\tau-t}\|^2
+c\int_{-\8}^0\me^{\lam_0\sig}\ol{\cQ}_1(\sig+\tau,\theta_{-\tau}\w)\di\sig\)\\
\leqslant&\me^{-\lam_0\vsig}\ol{\cR}^2(\lam_0,\tau,\w),
\end{align*}
which implies that $\ol{L}_2$ vanishes as $N\ra+\8$.
The term $\ol{L}_3$ surely converges to $0$ as $N\ra+\8$ by the limit \eqref{4.23}
and almost the same discussion in the second paragraph from bottom in the proof of
Lemma \ref{le3.5}.
Eventually, we have indeed proved the conclusion \eqref{4.25}.
\eo

Based on these estimates, we can hence obtain the existence of pullback random
attractors for the NRDS $\vp$ generated by \eqref{2.4n} with $k=0$.

\bt\label{th4.4m} Let the hypotheses \textbf{(F)}, \textbf{(G1)}, \textbf{(G2)},
\textbf{(H1)} and \textbf{(H2)} hold.
Then the NRDS $\vp$ defined in \eqref{4.15} has a unique pullback random $\ol{\cD}$-attractor
$\ol{\cA}$ and $\ol{\cA}\in\ol{\cD}$.
\et
\bo Similarly, we check the two conditions presented in Theorem \ref{th2.4n}.

Let $\ol{K}:=\{\ol{K}(\tau,\w)\}_{(\tau,\w)\in\R\X\W}$ be given in \eqref{4.19}.
We know that $\ol{K}$ is random by Proposition \ref{pro2.2n}.
Since $u_{\tau-t}\in D(\tau-t,\theta_{-t}\w)$ with $D\in\ol{\cD}$,
using \eqref{4.14} and \eqref{4.17},
we obtain that for each $\tau\in\R$ and $\bP$-a.s. $\w\in\W$,
\begin{align*}\|u(\tau,\tau-t,\theta_{-\tau}\w,u_{\tau-t})\|=&
\|v(\tau,\tau-t,\theta_{-\tau}\w,v_{\tau-t})+z(\w)h(\tau)\|\\
\leqslant&\ol{\cR}(\lam,\tau,\w)+|z(\w)|\|h(\tau)\|.\end{align*}
Therefore, $\ol{K}$ is a pullback $\ol{\cD}$-absorbing set for $\vp$ in $\ell^2$.
By \eqref{4.16} and \eqref{4.20}, we know $\ol{K}\in\ol{\cD}$.

We apply Lemma \ref{le4.3} to prove the pullback $\ol{\cD}$-asymptotic compactness of $\vp$.
For each $D\in\ol{\cD}$, consider $u_{\tau-t}\in D(\tau-t,\theta_{-t}\w)$.
By \eqref{4.15} and $v_{\tau-t}=\me^{-kz(\theta_{-t}\w)}(u_{\tau-t}+k^{-1}h(\tau-t))$,
we have
\begin{align*}&\sum_{|i|>N}\left|\(\vp(t,\tau-t,\theta_{-t}\w,u_{\tau-t})\)_i\right|^2
=\sum_{|i|>N}\left|v_i(\tau,\tau-t,\theta_{-\tau}\w,v_{\tau-t})+z(\w)h_i(\tau)\right|^2\\
\leqslant&2\(\sum_{|i|>N}\left|v_i(\tau,\tau-t,\theta_{-\tau}\w,v_{\tau-t})\right|^2
+z^{2}(\w)\sum_{|i|>N}|h_i(\tau)|^2\).
\end{align*}
On account of \eqref{4.25} and $h(\tau)\in\ell^2$, we obtain
the pullback $\cD$-asymptotic nullity of $\vp$ and
hence the pullback $\ol{\cD}$-asymptotic compactness.

The conclusion is thus proved by Theorem \ref{th2.4n},
and the proof is finished.
\eo
\subsection{Invariant sample measures and sample statistical solutions}

Similar to the treatment in Subsection \ref{ss3.4},
in order to attain the existence of invariant sample measures of $\vp$ in $\ell^2$
for the NRDS $\vp$ generated by \eqref{2.4n} with $k=0$,
we are to verify the boundedness of the mapping
$\tau\mapsto\vp(t-\tau,\tau,\theta_{\tau}\w,\phi)$ over $(-\8,t]$
and the joint continuity of the mapping $(\tau,\phi)\mapsto\vp(t-\tau,\tau,\theta_\tau\w,\phi)$
in $\ell^2$ for fixed $t\in\R$ and $\bP$-a.s. $\w\in\W$.

As before, we discuss this topic from $v(t,\tau,\w,\phi)$.
\bl\label{le4.5m} Let the hypotheses \textbf{(F)}, \textbf{(G1)}, \textbf{(G2)}, \textbf{(H1)}
and \textbf{(H2)} hold.
Then for each fixed $(t,\phi)\in\R\X\ell^2$, $\bP$-a.s., $\w\in\W$ and all $\tau\leqslant t$,
\be\label{4.27}\|v(t,\tau,\w,\phi)\|\leqslant\ol{\cM}(t,\w,\phi)<+\8,\ee
with
$$\ol{\cM}^2(t,\w,\phi):=\|\phi\|^2
+c\int_{-\8}^t\me^{\lam(\vsig-t)}\ol{\cQ}_1(\vsig,\w)\di\vsig.
$$
\el
\bo Choosing $\sig=\tau$, $s=t$ and $v_\tau=\phi$ in \eqref{4.11},
we have
\be\|v(t)\|^2\leqslant\me^{\lam(\tau-t)}\|\phi\|^2
+c\int_\tau^t\me^{\lam(\vsig-t)}\ol{\cQ}_1(\vsig,\w)\di\vsig,
\label{4.27m}\ee
which is uniformly dominated by $\ol{\cM}^2(t,\w,\phi)$ for each $\tau\leqslant t$.
The boundedness of $\ol{\cM}^2(t,\w,\phi)$ is obvious by repeating the estimates
for \eqref{4.18}.\eo

\bl\label{le4.6m} Let the hypotheses \textbf{(F)}, \textbf{(G1)}, \textbf{(G2)},
\textbf{(H1)} and \textbf{(H2)} hold.
Fix $(t,\phi)\in\R\X\ell^2$ and $\bP$-a.s., $\w\in\W$.
Then for every $\ve>0$, there exists $\ol{\de}_1=\ol{\de}_1(\ve,t,\w,\phi)>0$ such that when
$\tau\in(t-\ol{\de}_1,t)$,
$$\|v(t,\tau,\w,\phi)-\phi\|<\ve.$$
\el
\bo First, similar to \eqref{3.47z}, we also have
\be\label{4.28}\|v(t,\tau,\w,\phi)-\phi\|^2
\leqslant(t-\tau)\int_{\tau}^t\|\ol{F}(s,\w,v(s))\|^2\di s.\ee
Due to the following estimates:
\be\|\ol{F}_1(v)\|^2\leqslant c\|v\|^2,\ee
\begin{align}\|\ol{F}_2(s,\w,v(s))\|^2
\leqslant&c\|f(v(s)+z(\theta_s\w)h(s))\|^2+\|g(t,v(s)+z(\theta_s\w)h(s))\|^2\notag\\
\leqslant&c\(\|\al_1\|_1^{\frac{2q}{1+\epsilon}}+\|v+z(\theta_s\w)h(s)\|^{2q}
+\|\al_2\|_1^2\)\notag\\
&+c\(\|\psi_3(s)\|_\8^{2q}+\|v+z(\theta_s\w)h(s)\|^{2q}+\|\psi_4(s)\|^2\)\notag\\
\leqslant&c\(\|\al_1\|_1^{\frac{2q}{1+\epsilon}}+\|\al_2\|_1^2+\|\psi_3(s)\|_\8^{2q}
+\|\psi_4(s)\|^2\)\notag\\
&+c\|v(s)\|^{2q}+c|z(\theta_s\w)|^{2q}\|h(s)\|^{2q},
\end{align}
\be\|\ol{F}_3(s,\w)\|^2\leqslant cz^2(\theta_s\w)\(\|h(s)\|^2+\|h'(s)\|^2\)\ee
and by \eqref{4.27}, for $s\in[\tau,t]$,
\be\|v(s)\|^2\leqslant\ol{\cM}^2(s,\w,\phi)
=\|\phi\|^2+c\int_{-\8}^s\me^{\lam(\vsig-s)}\ol{\cQ}_1(\vsig,\w)\di\vsig
\leqslant\me^{\lam(t-s)}\ol{\cM}^2(t,\w,\phi),\label{4.32}
\ee
and by setting
\begin{align*}\ol{\cQ}_2(s,\w):=&\|\al_1\|_1^{\frac{2q}{1+\epsilon}}+\|\al_2\|_1^2
+\|\psi_3(s)\|_\8^{2q}+\|\psi_4(s)\|^2+|z(\theta_s\w)|^{2q}\|h(s)\|^{2q}\\
&+z^2(\theta_s\w)\(\|h(s)\|^2+\|h'(s)\|^2\),\end{align*}
we deduce by the inequalities from \eqref{4.28}-\eqref{4.32} that when $t>\tau>t-1$,
\begin{align}\frac{1}{t-\tau}\|v(t,\tau,\w,\phi)-\phi\|^2
\leqslant&c\int_{\tau}^t\(\|v(s)\|^2+\|v(s)\|^{2q}+\ol{\cQ}_2(s,\w)\)\di s\notag\\
\leqslant&c\ol{\cM}^2(t,\w,\phi)\int_{t-1}^t\me^{\lam(t-s)}\di s
+c\ol{\cM}^{2q}(t,\w,\phi)\int_{t-1}^t\me^{\lam q(t-s)}\di s\notag\\
&+c\int_{t-1}^t\ol{\cQ}_2(s,\w)\di s.
\label{4.34}\end{align}
Similarly, one can easily deduce the boundedness of the right side of \eqref{4.34}
by the continuity of $h$.

Consequently, for all $\ve>0$, we have a sufficiently small
$\ol{\de}_1:=\ol{\de}_1(\ve,t,\w,\phi)\in(0,1)$ such that when $\tau\in(t-\ol{\de}_1,t)$,
$\|v(t,\tau,\w,\phi)-\phi\|^2<\ve^2$. The proof is then finished.
\eo
\bl\label{le4.7m} Let the hypotheses \textbf{(F)}, \textbf{(G1)}, \textbf{(G2)},
\textbf{(H1)} and \textbf{(H2)} hold.
Then given each fixed $\tau,\,t\in\R$ with $t>\tau$, $\phi\in\ell^2$
and $\bP$-a.s. $\w\in\W$, for every $\ve>0$,
there exists $\ol{\de}_2=\ol{\de}_2(\ve,t,\tau,\w,\phi)>0$
such that whenever $|\tau'-\tau|+\|\phi'-\phi\|<\ol{\de}_2$,
$$\|v(t,\tau',\w,\phi')-v(t,\tau,\w,\phi)\|<\ve.$$
As a result, the mapping $(\tau,\phi)\mapsto v(t,\tau,\w,\phi)$ is continuous
for each $t\geqslant\tau$ and $\bP$-a.s. $\w\in\W$.
\el

\bo The proof is parallel to that of Lemma \ref{le3.9} and so we omit it here.
\eo

\bl\label{le4.8m} Let the hypotheses \textbf{(F)}, \textbf{(G1)}, \textbf{(G2)},
\textbf{(H1)} and \textbf{(H2)} hold.
Then the mapping $\tau\!\mapsto\! u(t,\tau,\w,\phi)$ is bounded for $\tau\in(-\8,t]$,
and the mapping $(\tau,\phi)\ra u(t,\tau,\w,\phi)$ is continuous.\el
\bo By \eqref{4.14} and Lemma \ref{le4.5m},
we know that for all $\tau\in(-\8,t]$,
\begin{align*}&\|u(t,\tau,\w,\phi)\|=
\|v(t,\tau,\w,\phi-z(\theta_\tau\w)h(\tau))+z(\theta_t\w)h(t)\|\\
\leqslant&\|v(t,\tau,\w,\phi-z(\theta_\tau\w)h(\tau))\|
+|z(\theta_t\w)|\|h(t)\|.\end{align*}
By \eqref{4.27m}, we infer that
\begin{align*}&\|v(t,\tau,\w,\phi-z(\theta_\tau\w)h(\tau))\|^2\\
\leqslant&\me^{\lam(\tau-t)}\|\phi-z(\theta_\tau\w)h(\tau)\|^2
+c\int_\tau^t\me^{\lam(\vsig-t)}\ol{\cQ}_1(\vsig,\w)\di\vsig\\
\leqslant&c\me^{\lam(\tau-t)}\|\phi\|^2+c\me^{\lam(\tau-t)}\|h(\tau)\|^2
+c\int_\tau^t\me^{\lam(\vsig-t)}\ol{\cQ}_1(\vsig,\w)\di\vsig,
\end{align*}
which is uniformly bounded for all $\tau\in(-\8,t]$ by \eqref{3.27z},
\eqref{4.23} and the hypotheses \textbf{(G2)} and \textbf{(H2)}.

The continuity of the mapping $(\tau,\phi)\ra u(t,\tau,\w,\phi)$ follows
from the continuity of $z(\theta_\tau\w)$ and $h(\tau)$ in $\tau$
and Lemma \ref{le4.7m}. The proof is complete.\eo

Based on Proposition \ref{pro2.6n}, Theorem \ref{th4.4m}, Lemma \ref{le4.8m}
and \eqref{4.15},
we can construct the invariant sample measures of $\vp$ on $\ell^2$ as follows.

\bt\label{th4.9m}
Suppose that the hypotheses \textbf{(F)}, \textbf{(G1)}, \textbf{(G2)}, \textbf{(H1)}
and \textbf{(H2)} hold.
Let $\ol{\cA}$ be the unique pullback random $\ol{\cD}$-attractor of $\vp$
given in Theorem \ref{th4.4m}.
Then for a given generalized Banach limit $\disp\LIM_{t\ra-\8}$ and
a mapping $\xi:\R\X\W\ra\ell^2$ with $\{\xi(\tau,\w)\}_{(\tau,\w)\in\R\X\W}\in\ol{\cD}$
such that $\xi(\tau,\theta_\tau\w)$ is continuous in $\tau$,
there exists a family of Borel probability measures $\{\ol{\mu}_{\tau,\w}\}_{(\tau,\w)\in\R\X\W}$
on $\ell^2$ such that the support of the measure $\ol{\mu}_{\tau,\w}$ is contained
in $\ol{\cA}(\tau,\w)$ and for all $\Upsilon\in\cC(\ell^2)$,
\begin{align}&
 \LIM_{\tau\ra-\8}\frac{1}{t-\tau}\int_{\tau}^{t}
 \Upsilon(\vp(t-s,s,\theta_s\w,\xi(s,\theta_s\w)))\di s
 =\int_{\ol{\cA}(t,\theta_t\w)}\Upsilon(u)\ol{\mu}_{t,\theta_t\w}(\di u)\notag\\
=&\int_{\ell^2}\Upsilon(u)\ol{\mu}_{t,\theta_t\w}(\di u)
 =\LIM_{\tau\ra-\8}\frac{1}{t-\tau}\int_{\tau}^{t}\int_{\ell^2}\Upsilon(\vp(t-s,s,\theta_s\w,u))
 \ol{\mu}_{s,\theta_s\w}(\di u)\di s.
\label{4.35}\end{align}
Additionally, $\ol{\mu}_{\tau,\w}$ is invariant in the sense that for all $\tau\in\R$ and
$\bP$-a.s. $\w\in\W$,
\be\label{4.36}
 \int_{\ol{\cA}(\tau+t,\theta_t\w)}\Upsilon(u)\ol{\mu}_{\tau+t,\theta_t\w}(\di u)=
 \int_{\ol{\cA}(\tau,\w)}\Upsilon(\vp(t,\tau,\w,u))\ol{\mu}_{\tau,\w}(\di u),
 \hs\mb{for all }t\geqslant0.
\ee\et

At last, with the almost same argument for the case when $k>0$,
we can similarly conclude that the invariant sample measures constructed in Theorem \ref{th4.9m}
satisfies the stochastic Liouville type equation \eqref{1.7} and is a sample statistical solution,
which is stated specifically as follows.
\bt Suppose that the hypotheses \textbf{(F)}, \textbf{(G1)}, \textbf{(G2)},
\textbf{(G3)}, \textbf{(H1)} and \textbf{(H2)} hold.
Let $V$, $H$ and $V^*$ be defined as \eqref{3.59z1} and $\tilde{F}$ and $\tilde{G}$
as \eqref{3.59z2} with $k=0$.
Let $\ol{\cA}$ be the pullback random $\ol{\cD}$-attractor of $\vp$ given in Theorem \ref{th4.4m}
and $\{\ol{\mu}_{\tau,\w}\}_{(\tau,\w)\in\R\X\W}$ be the invariant sample measures constructed in
Theorem \ref{th4.9m}.
Then the invariant sample measures $\{\ol{\mu}_{\tau,\w}\}_{(\tau,\w)\in\R\X\W}$ satisfy
the stochastic Liouville type equation \eqref{1.7} for each $\Psi\in\cT$
with $\cT$ given in Definition \ref{de1.1}.

Moreover, the family of Borel probability measures
$\{\ol\mu_{\tau,\w}\}_{(\tau,\w)\in\R\X\W}$ obtained
in Theorem \ref{th4.9m} is a sample statistical solution
of the equation \eqref{2.4n}.\et

\section{Summary and Remarks}

In this work, based on the existence of invariant sample measures,
we have actually extended the stochastic Liouville type equation
established by Chen and Yang in \cite{CY23} to a more general framework,
and further given a suitable definition of sample statistical solution with
respect to the invariant sample measures for general stochastic differential equations.

As an application, we construct the invariant sample measures
for the nonautonomous stochastic lattice Cahn-Hilliard equation with
nonlinear noise and validate that the invariant sample measures is a sample statistical solution
in the abstract framework.
In this situation, the necessary hypotheses endowed upon $f$ and $g$ lead to
the equivalence relation between $H$ and $V$ (set in Section \ref{s1}),
which thereby helps largely to guarantee related integrability.

For further applications, it is a more interesting theme to establish sample statistical solutions
for more general situations, such as lattice systems on weighted spaces,
partial differential equations on bounded smooth domains and unbounded domains.

\section*{Data Availability}

Data sharing is not applicable to this article as no new data were created
or analyzed in this study.

\section*{Acknowledgements}

Our work was supported by grants from the National Natural Science Foundation of China
(NNSFC Nos. 11801190 and 12101462).

\end{document}